\input ./style/arxiv-general.cfg
\documentclass[aos,MSNbibl,dvips]{arximspdf}
\makeatletter
   \@ifpackageloaded{graphicx}{}{\usepackage{graphicx}}
\makeatother
\usepackage{algorithm}
\usepackage{mathrsfs}

\doi{10.1214/15-AOS1346}% Updated by VTEXPTS2LaTeX.exe, 02.07.2015
\volume{43}
\issue{5}
\pubyear{2015}
\firstpage{2296}
\lastpage{2325}
\docsubty{FLA}

\makeatletter
\def\sfrac#1#2{#1/#2}

\newcommand{\rrvert}{\vert}
\newcommand{\rrVert}{\Vert}
\newcommand{\llvert}{\vert}
\newcommand{\llVert}{\Vert}
\newcommand{\overset}{\stackrel}
\newproclaim{assumption}{Condition}%[section]
\newtheorem{proposition}{Proposition}%[section]
\newtheorem{theorem}{Theorem}
\newtheorem{corollary}{Corollary}
\newtheorem{lemma}{Lemma}
\newcommand{\btheta}{\bolds\theta}
\newcommand{\bbeta}{{\bolds\beta}}
\newcommand{\bfeta}{\bolds\eta}
\newcommand{\bxi}{\bolds\xi}
\newcommand{\bY}{{\mathbf Y}}
\newcommand{\bb}{\mathbf{b}}
\newcommand{\bR}{{\mathbf R}}
\newcommand{\bX}{{\mathbf X}}
\newcommand{\bc}{{\mathbf c}}
\newcommand{\bv}{{\mathbf v}}
\newcommand{\bh}{{\mathbf h}}
\newcommand{\bP}{{\mathbf P}}
\newcommand{\F}{{\mathcal F}}
\def\RR{\mathbb{R}}
\newcommand{\bff}{{\mathbf f}}
\newcommand{\hf}{\hat{\bff}}
\newcommand{\hbfeta}{\hat{\bfeta}}
\newcommand{\bzero}{\mathbf{0}}
\newcommand{\Mfrak}{\mathfrak{M}}
\newcommand{\hbff}{\hat{\bff}}
\newcommand{\hbxi}{\hat{\bxi}}
\newcommand{\veps}{\varepsilon}
\def\toD{\overset{\mathscr{D}}{\longrightarrow}}
\makeatother

\begin{document}
\begin{frontmatter}

\title{Functional additive regression}
\runtitle{Functional additive regression}

\begin{aug}
% Corresponding author: Yingying Fan - fanyingy@usc.edu% Updated by
%VTEXPTS2LaTeX.exe, 02.07.2015 11:31
\author[A]{\fnms{Yingying}~\snm{Fan}\corref{}\thanksref{T1}\ead[label=e1]{fanyingy@marshall.usc.edu}},
\author[A]{\fnms{Gareth M.}~\snm{James}\ead[label=e2]{gareth@marshall.usc.edu}}
\and
\author[A]{\fnms{Peter}~\snm{Radchenko}\thanksref{T2}\ead[label=e3]{radchenk@marshall.usc.edu}}
\runauthor{Y. Fan, G.~M. James and P. Radchenko}
\affiliation{University of Southern California}
%\dedicated{}
\address[A]{Data Sciences and Operations Department\\
Marshall School of Business\\
University of Southern California\\
Los Angeles, California 90089\\
USA\\
\printead{e1}\\
\phantom{E-mail: }\printead*{e2}\\
\phantom{E-mail: }\printead*{e3}}
\end{aug}
\thankstext{T1}{Supported by NSF CAREER Award DMS-11-50318 and USC
Marshall summer research funding.}
\thankstext{T2}{Supported in part by NSF Grant DMS-12-09057.}

% HISTORY:
%
\received{\smonth{11} \syear{2014}}% Updated by VTEXPTS2LaTeX.exe,
%02.07.2015 11:31
%
\revised{\smonth{5} \syear{2015}}% Updated by VTEXPTS2LaTeX.exe,
%02.07.2015 11:31

% ABSTRACT
%
\begin{abstract}
We suggest a new method, called \textit{Functional Additive
Regression}, or
FAR, for efficiently performing high-dimensional functional regression.
FAR extends the usual linear regression model involving a functional
predictor, $X(t)$, and a scalar response, $Y$, in two key respects.
First, FAR uses a penalized least squares optimization approach to
efficiently deal with high-dimensional problems involving a large number
of functional predictors. Second, FAR extends beyond the
standard linear regression setting to fit general nonlinear additive
models. We demonstrate that FAR can be implemented with a wide range of
penalty functions using a highly efficient coordinate descent
algorithm. Theoretical results are developed which provide motivation
for the FAR optimization criterion. Finally, we show through simulations
and two real data sets that FAR can significantly outperform competing methods.
\end{abstract}

% KEYWORDS
\end{frontmatter}

%s1 #&#
\section{Introduction}
The univariate functional regression situation, where one models the
relationship between a scalar response, $Y$, and a functional predictor,
$X(t)$, has recently received a great deal of attention. A few examples
include \cite{hastie7,hall1,alter1,hall2,james3,cardot1,ferraty2,james5,muller1,ferraty3,chen2}. See Chapter~15 of \cite{Bramsay3} for a thorough
discussion of the issues involved with fitting such data.

Most
work in this area involves different approaches for fitting the functional
linear regression model,
%
%e1 #&#
\begin{equation}
\label{flin} Y_i = \int\beta(t) X_i(t) \,dt +
\varepsilon_i, \qquad i = 1,\ldots, n.
\end{equation}
For notational convenience, we assume throughout this paper that the
response and predictors have been centered so the intercept can be
ignored. Model (\ref{flin}) provides a natural extension of linear
regression to
the functional domain but it has two significant limitations. First,
it assumes a single predictor, while functional regression situations
involving a large number of predictors, $X_{i1}(t), X_{i2}(t), \ldots,
X_{ip}(t)$, are becoming increasingly common. For example, \cite{storey1}
analyzes two gene expression data sets measured over time, which
involve only a small number of patients but tens of thousands of
functional predictors. Second, (\ref{flin}) is relatively inflexible
because it assumes a
linear relationship between the predictor and response. Just as in the
standard regression setting more accurate fits can often be produced by
modeling a nonlinear relationship.

In this paper, we address both of these limitations using a functional
additive regression framework of the form
%
%e2 #&#
\begin{equation}
\label{far} Y_i = \sum_{j=1}^p
f_j (X_{ij} )+ \varepsilon_i, \qquad i = 1,
\ldots, n,
\end{equation}
where the $f_j$'s are general nonlinear functions of $X_{ij}(t)$.
There has been some previous work extending the classical functional
regression model. James and Silverman \cite{james5} proposed an index model to
implement a nonlinear functional \mbox{regression}, and, more recently, both
\cite{ferraty4} and \cite{chen2} extended this work to a fully
nonparametric setting and provided further theoretical motivation.
However, all of these approaches are primarily intended for the
univariate setting, where $p=1$. Lian \cite{lian1} did consider a multivariate
setting involving both functional and scalar predictors, but with only
a single functional predictor, so the corresponding model does not
extend to (\ref{far}).
%\cite{zhu1} propose a functional variable selection approach using a
%constrained least squares method with similarities to the group lasso,
%while \cite{zhu2} utilizes a Bayesian formulation. Both methods
%provide powerful extensions to the generalized linear models domain
%but do not consider nonlinear relationships.
James and Silverman \cite{james5} proposed a kernel based method for fitting~(\ref{far}), which works well in low-dimensional situations. However,
they do not attempt to perform any kind of variable selection. As a
result, the method suffers from computational and statistical issues
when $p$ is large, such as for the gene expression data in \cite
{storey1}. Zhu et al. \cite{zhu2} proposed a Bayesian variable selection
approach for selecting and estimating important functional predictors
in a classification setting.
%, where the functional predictors were modeled using orthogonal basis
%functions.
However, while their method can potentially be implemented on a large
number of functions, it still assumes a linear relationship between the
response and predictors. Finally, a recent paper \cite{fan1} considers
a more general form of (\ref{far}) where the response is also
functional. Their approach appears to work well but the paper does not
provide any theoretical results. See also \cite
{Febrero-Bande,Goia12,Mas} for additional recent developments on
functional regression
models with multiple functional covariates under various model settings.

Fitting
(\ref{far}) in the high-dimensional setting poses
a couple of significant complications. First, in order to make the problem
feasible, we must assume sparsity in the predictor space, that is, that most
of the predictors are unrelated to the response.
% so $\left\Vert f_j\right\Vert=0$ for most
%$j$. In this paper we take $\left\Vert f_j\right\Vert$ as the $L_2$
%norm of $f_j$ over
%the
%distribution of $X_j$ i.e. $\left\Vert f_j\right\Vert^2=Ef^2_j(X_j)$.
Thus, we need an
approach that can automatically perform high-dimensional variable
selection on nonlinear functions. Second, (\ref{far}) involves
estimating functions, $f_j(x)$, of functional predictors, $X_{ij}(t)$.
Even in the
univariate situation, involving a single predictor, there has been little
research on this problem and the best approach is unclear. Most current methods
involve using the first few functional principal component scores of
$X_{ij}(t)$ as a finite-dimensional predictor space \cite{muller2}.
However, the principal component scores are computed independently from
the response, in an unsupervised fashion, so there is no a priori
reason to believe that these scores
will correspond to the best dimensions for the regression problem.

In this paper, we suggest a new penalized least squares method called
\emph{Functional Additive Regression}, or FAR, for fitting a
nonlinear functional additive model. FAR makes three important
contributions. First, it efficiently fits high-dimensional functional
models while simultaneously
performing variable selection to identify the relevant predictors. This
is an area
that has historically received very little attention in the functional
domain, but the importance of the connections between functional and
high-dimensional statistics are just starting to become clear. See, for
example, the recent conference on this topic \cite{Bongiorno1}.

Second,
FAR extends beyond the standard linear regression setting to fit general
nonlinear additive models.
%In order to fit (\ref{far}) one must place some form of restriction on
%the $f_j$'s.
FAR models $f_j(x)$ as a nonlinear function of a one-dimensional
linear projection of $X_{ij}(t)$; a functional version of the \emph
{single index model} approach. Our method uses a supervised fit to
automatically project the functional predictors into the best
one-dimensional space. We believe this is an important distinction
because projecting into the unsupervised PCA space is currently the
dominant approach in functional regressions, even
though it is well known that this space need not be optimal for
predicting the response.

Third, FAR can be implemented using a wide range of penalty functions and
a highly efficient coordinate descent algorithm. In the linear case, we
establish a number of theoretical results, which show that, under suitable
conditions and for an appropriately chosen penalty function, FAR is
guaranteed to asymptotically choose the correct model as $n$ and $p$ go to
infinity. Theoretical investigation for the nonlinear FAR approach
presents some serious additional challenges, because the regression
functions, $f_j$, are estimated rather than known. We allow the number of
functional predictors, $p$, to grow faster than the number of
observations, $n$, and establish asymptotic bounds on the $\ell_2$
estimation error for each of the estimated regression functions. The
difficulties associated with the high-dimensional nature of the functional
data are exacerbated by the large number of estimated components in the
additive regression model for the response. Moreover, the functional
aspect of the data (infinite dimensional predictors) adds further
complexity to the already very challenging problem. Our method of proof
uses ideas from the estimation theory for high-dimensional additive models
\cite{Bbuhlmann1,huang1,meier1}. However, the proof itself is new, rather than a compilation
of existing results.

%Contributions:
%\begin{enumerate}
%\item Variable selection with $p$ functional predictors where $p$ is
% large. Very little work with multivariate functional data.
%\item Functional regression with additive nonlinear response
%surfaces. A
% natural extension of GAM to functional predictors.
%\item Methodology and theory for arbitrary penalty functions.
%\end{enumerate}

Our paper is set out as follows. In Section~\ref{linFARsec}, we
develop the FAR method for performing high-dimensional functional
regression. Section~\ref{funcindecsec} uses functional index models
to motivate the FAR model. Then Section~\ref{linearfarsec} presents
the optimization criterion and an efficient coordinate
descent algorithm for fitting FAR in the linear regression setting.
Finally, Section~\ref{nlfar} extends the algorithm to the nonlinear
regression framework.
% and is presented which can be solved using
%Our main theoretical results are provided in Section~\ref{theorysec}.
In Section~\ref{theorysec}, we provide a number of theoretical
results. We first prove that, under appropriate conditions, the linear
version of FAR will asymptotically include all the true signal
variables and remove all the noise
predictors from the model. In addition, we provide an asymptotic bound
on the estimation error of the signal functions, $f_j(x)$, under the
vector infinity norm,\vadjust{\goodbreak} and show that the FAR estimator is asymptotically
normal. In the nonlinear setting, we establish the rate of
convergence, with respect to the $\ell_2$ distance, for the estimates
of the regression functions, $f_j(x)$, corresponding to each of the predictors.
We also investigate the variable selection properties of our estimator
and show that, under some conditions, it can recover the index set of
the signal predictors.
%Section~\ref{nonlinfarsec} extends FAR to the general nonlinear
%additive model framework using a
%supervised approach for projecting the predictors into a finite
%dimensional space.
Extensive simulation results are presented in Section~\ref{simsec}. We
compare FAR to other functional regression methods and demonstrate its
superior performance in many settings. Finally, we apply FAR to both
medium and high-dimensional real data sets in Section~\ref{realsec},
and end with a discussion in
Section~\ref{discsec}.
%Proofs of all theoretical results are provided in a separate online
%appendix.

%s2 #&#
\section{Functional additive regression}
\label{linFARsec}

Let $\bff_j= (f_j(X_{1j}),\ldots,f_j(X_{nj}) )^T$. Then
our general approach for fitting (\ref{far}) is to minimize the
following penalized regression criterion over $\bff_1, \bff_2,\ldots,
\bff_p$:
%
%e3 #&#
\begin{equation}
\label{fundamental} \frac{1}{2n} \Biggl\llVert\bY-\sum
_{j=1}^p \bff_j \Biggr\rrVert
^2_2 + \sum_{j=1}^p
\rho_{\lambda_n} \biggl(\frac{1}{\sqrt{n}}{\llVert\bff_j\rrVert
}_2 \biggr),
\end{equation}
where $\bY=(Y_1,\ldots, Y_n)^T$, $\rho_{\lambda_n}(t)$ is a penalty
function, $\lambda_n$ is the regularization parameter and ${\llVert
\bff_j\rrVert
}_2=\sqrt{\bff_j^T\bff_j}$. To aid the presentation, we drop the
subscript and use $\llVert \cdot\rrVert $ to denote the $\ell_2$ norm
of a vector
in the future. Although it may not be immediately obvious from this
formulation, we show that minimizing (\ref{fundamental}) will in
general automatically implement variable selection by shrinking a
subset of the $\bff_j$'s to exactly zero. In this article, we explore
general concave functions for $\rho$, with the $\ell_1$ penalty $\rho
_{\lambda}(t) = \lambda t$ considered as a special case. There is by
now a substantial literature demonstrating the
%It has been justified both theoretically and empirically by many
%researchers that concave
%penalty functions have
advantages of concave penalty functions for high-dimensional problems
\cite{fan2,fan3,lv1,fan5,loh13}.

We assume that the trajectories of functional predictors, $X_{ij}(t)$,
are fully observed. Our methodology and theoretical results can be
extended to the case of densely observed predictors under additional
smoothness and regularity assumptions. However, for the clarity of the
exposition we do not investigate this case in the paper.

%\subsection{Linear FAR Criterion}

%s2.1 #&#
\subsection{Functional index models}
\label{funcindecsec}

Minimizing (\ref{fundamental}) requires specifying the form of
$f_j(x)$. A limitation of linear functional regression models is that
they can
perform poorly when there is a nonlinear relationship between $X(t)$ and
$Y$. However, the infinite-dimensional nature of $X(t)$ makes it
challenging to model a nonlinear relationship between the predictor and
response. As a result, relatively few papers have investigated
this extension. Most methods focus on approximating $X(t)$ using its first
few functional principal components and then implementing nonlinear fits
using the principal component scores as predictors \cite{muller2}. However,
this unsupervised approach has the usual limitation; the directions which
explain $X(t)$ best may not be the most appropriate for predicting the
response.

In the multivariate setting, index models are commonly used for
providing nonlinear fits to high-dimensional data. For a centered
response, the standard single index model can be expressed in the form
$Y=g(\bbeta^T\bX) + \varepsilon$,
where $g(x)$ is a general
nonlinear function and $\bbeta$ is a norm one vector representing the
best single direction to project the predictors into. A key advantage of
the index model formulation is that $\bbeta$ is chosen in a supervised
fashion, incorporating both the response and predictors, potentially
providing more accurate fits. Index models can be naturally extended to
functional predictors using the formulation $f_j(X_{ij}) = g_j
(\int\beta_j(t)X_{ij}(t)\,dt )$, where $g_j(x)$ and $\beta_j(t)$
are both nonparametric smooth functions, and the integral is
well-defined. %Since $\beta_j(t)$ corresponds to a direction that we
%project $X_{ij}(t)$ into we impose the constraint $\left\Vert\beta_j
%\right\Vert=1$.
Functional single index models have been considered previously. For
example, \cite{james5,amato1,ait1,chen2,ferraty4}, all fit index models to functional data, but these
previous approaches all concentrate on the $p=1$ problem.

Using this nonlinear representation, the FAR model (\ref{far}) can be
expressed as
%
%e4 #&#
\begin{equation}
\label{simfar} Y_i = \sum_{j=1}^p
g_j \biggl(\int\beta_j(t)X_{ij}(t)\,dt \biggr)
+ \varepsilon_i.
\end{equation}
For identifiability, in addition to centering the response, we also
center the regression functions: $\sum_{i=1}^n g_j (\int\beta
_j(t)X_{ij}(t)\,dt )=0$ for all~$j$. Note that index functions
$\beta_j$ are only identifiable up to multiplications by nonzero
constants, however, our focus is on estimating $f_j$ rather than $\beta
_j$. The general FAR optimization criterion~(\ref{fundamental}) becomes
%
%e5 #&#
\begin{equation}
\label{nonlinfar} \frac{1}{2n} \Biggl\llVert\bY-\sum
_{j=1}^p g_j \biggl(\int
\beta_j(t) \bX_j(t)\,dt \biggr)\Biggr\rrVert
^2 + \sum_{j=1}^p
\rho_{\lambda_n} \biggl(\frac{1}{\sqrt n}\llVert\bff_j\rrVert
\biggr),
\end{equation}
where\vspace*{1pt} $\bX_j(t) = (X_{1j}(t), \ldots, X_{nj}(t))^T$ and $g_j
(\int\beta_j(t) \bX_j(t)\,dt ) = (f_j(X_{1j}), \ldots,\break  f_j(X_{nj}))^T$.

%s2.2 #&#
\subsection{Linear FAR}
\label{linearfarsec}
Our approach for minimizing (\ref{nonlinfar}) is easiest to understand
by first considering the situation where $f_j(x)$ is taken to be
linear. Hence, in this section we develop FAR in the setting where
$g_j(x)$ is set to the identity function, in which case FAR reduces to
a multivariate functional linear regression model.
%We further extend FAR to the more general nonlinear setting in
%Section~\ref{nlfar}.

%s2.2.1 #&#
\subsubsection{FAR criterion}

We assume without loss of generality that each predictor is observed
over the range $0\le t \le1$. Hence, in the linear setting,
%
%e6 #&#
\begin{equation}
\label{linearf} f_j(X_{ij}) = \int_0^1
\beta_j(t)X_{ij}(t)\,dt,
\end{equation}
where $\beta_j(t)$ is an unknown smooth coefficient function, and the
FAR optimization criterion becomes
%
%e7 #&#
\begin{equation}
\label{linearfar} \frac{1}{2n} \Biggl\llVert\bY-\sum
_{j=1}^p \int_0^1
\beta_j(t) \bX_j(t)\,dt\Biggr\rrVert^2 +
\sum_{j=1}^p \rho_{\lambda_n} \biggl(
\frac{1}{\sqrt{n}}\llVert\bff_j\rrVert\biggr),
\end{equation}
where $\bX_j(t)= (X_{1j}(t),\ldots,X_{nj}(t) )^T$.

%Specifically we
%select a $q$-dimensional orthonormal basis function, $\bb(t)$,
%satisfying $\int_0^1\bb(t)\bb^T(t)\,dt = I_q$ with $I_q$ the $q\times q$
%identity matrix. Hence the coefficient
%functions and predictors can be expressed as
Given an orthonormal basis $\{b_{l}(t)\}$, the functional predictors
and the corresponding regression coefficients can be decomposed as
%
%e8 #&#
\begin{equation}
\label{eqbasis-rep} X_{ij}(t) = \sum_{l=1}^\infty
\theta_{ijl}b_l(t), \qquad\beta_j(t) = \sum
_{l=1}^\infty\eta_{0,jl}b_l(t),
\end{equation}
where $\theta_{ijl}$ and $\eta_{0,jl}$ are the coefficients of
$X_{ij}(t)$ and $\beta_j(t)$ corresponding to the $l$th basis function
$b_l(t)$, respectively. Using (\ref{eqbasis-rep}), the $j$th additive
component has the following representation
%
%e9 #&#
\begin{equation}
\label{fexpan} f_j(X_{ij}) = \int_0^1
X_{ij}(t)\beta_j(t)\,dt = \sum_{l=1}^\infty
\theta_{ijl}\eta_{0,jl}.
\end{equation}

In order for the functions optimizing (\ref{linearfar}) to have
nontrivial solutions, some form of smoothness constraint must be
imposed on
the $\beta_j(t)$'s. Two standard approaches are to include a smoothness
penalty in the optimization criterion or alternatively to restrict the
functions to some low-dimensional class. In this setting, either
approach could be adopted but we use the latter method.
Specifically, for a given sequence of integers $q_n = o(n)$ depending
only on the sample size $n$, write $\bfeta_{0j} = (\eta_{0,j1},\ldots
,\eta_{0,jq_n})^T$ and $\btheta_{ij} = (\theta_{ij1},\ldots, \theta
_{ijq_n})^T$. Thus, the $j$th additive component $f_j(X_{ij})$ can be
approximately as $\btheta_{ij}^T\bfeta_{0j}$. Denote by $e_{ij}$ the
approximation error, that is,
%
%e10 #&#
\begin{equation}
\label{deferr} e_{ij} = f_j(X_{ij})-
\btheta_{ij}^T\bfeta_{0j} = \sum
_{l=q_n+1}^\infty\theta_{ijl}\eta_{0,jl}.
\end{equation}
Then by the Cauchy--Schwarz inequality and Condition~\ref{as3} in
Appendix~\ref{appb}, uniformly across all $i = 1,\ldots, n$ and $j\in
\Mfrak_0$,
%
%e11 #&#
\begin{eqnarray}
\qquad\llvert e_{ij}\rrvert^2& \leq&\sum
_{l=q_n+1}^\infty\eta_{0,jl}^2l^{-4}
\sum_{l=q_n+1}^\infty\theta_{ijl}^2l^{4}
\leq C^2 q_n^{-4}\sum
_{l=q_n+1}^\infty\eta_{0,jl}^2 \leq
\widetilde C C^2q_n^{-4}, \label{eqapp-err}
\end{eqnarray}
where $C$ and $\widetilde C$ are two positive constants defined in
Condition~\ref{as3}. Thus, for large enough $q_n$, the approximation
error is uniformly small.

%\begin{equation}
%\label{basisrep}
%\beta_j(t)=\bb(t)^T\bfeta_j + \widetilde r_j(t)\qquad
%\text{and}\quad X_{ij}(t)=\bb(t)^T\btheta_{ij} + r_{ij}(t),
%\end{equation}
%where $\widetilde r_j(t)$ and $r_{ij}(t)$ are remainder terms representing
%the bias in our basis approximation.
%The assumption that $\beta_j(t)$ and $X_{ij}(t)$ can be expressed
%using the same basis is made for simplicity of exposition. All
%the FAR calculations can be extended to the situation where the bases
%differ at the cost of some additional notation.

%where $e_{ij} = \int_0^1 r_{ij}(t)\beta_j(t)\,dt + \int_0^1\widetilde
%r_j(t)X_{ij}(t)\,dt - \int_0^1 r_{ij}(t)\widetilde r_j(t)\,dt$.
Let $\Theta_j$ be an
$n\times q_n$ matrix whose rows are formed by $\{\btheta_{ij},
i=1,\ldots, n\}$. Then, if $q_n$ is large enough, $f_j(X_{ij})\approx
\btheta_{ij}^T\bfeta_{0j}$ and (\ref{linearfar}) can be approximated by
%
%e12 #&#
\begin{equation}
\label{newopt2} \frac{1}{2n}\Biggl\llVert\bY- \sum
_{j=1}^p\Theta_j\bfeta_j
\Biggr\rrVert^2+\sum_{j=1}^p
\rho_{\lambda_n} \biggl(\frac{1}{\sqrt{n}}\llVert\Theta_j
\bfeta_j\rrVert\biggr).
\end{equation}
Note that the $\bfeta_j$'s must be estimated, but the
$\Theta_{j}$'s are calculated from the fully observed trajectories of
the functional predictors, $X_{ij}(t)$. Hence, we fit FAR by minimizing
(\ref{newopt2}) over $\bfeta_1,\ldots, \bfeta_p$.

%s2.2.2 #&#
\subsubsection{FAR algorithm}\label{seclinearfar}

The criterion given by (\ref{newopt2}) is still $p\times q_n$
dimensional, so is potentially challenging to optimize over, even if
$p$ is only of moderate size. However, in this form our FAR criterion
is closely related to the standardized group lasso \cite{simon1}
which allows us to develop an efficient algorithm to fit FAR.
In particular, a distinct advantage of (\ref{newopt2}) is that, when
using the Lasso penalty $\rho_{\lambda_n}(t)=\lambda_nt$, there is a
simple closed form expression for computing its minimum over $\bfeta_j$.

%pr1 #&#
\begin{proposition}\label{propL1}
If $\rho_{\lambda_n}(t) = \lambda_nt$, then the solution to (\ref
{newopt2}) satisfies $\hat{\bff}_j = \Theta_j\hbfeta_j$
%$$\hat{\bff}_j = \Theta_j\hat{\bfeta}_j= %$$
where
\[
\hbfeta_j = \biggl(1-\frac{\sqrt{n}\lambda_n}{\llVert S_j\bR_j\rrVert
} \biggr)_+\bigl(
\Theta_j^T \Theta_j\bigr)^{-1}
\Theta_j^T\bR_j,
\]
$S_j = \Theta_j(\Theta_j^T\Theta_j)^{-1}\Theta_j^T$, $\bR_j = \bY
- \sum_{k\neq j}\Theta_k\hbfeta_k$, and $z_+ = \max(0,z)$
represents the positive part of $z$.
\end{proposition}

The derivation of Proposition~\ref{propL1} involves simple algebra
and similar results are proved in \cite{ravikumar1} and \cite{simon1}
so we do not provide the proof here.
%
%Recent research has shown that coordinate descent algorithms can be an
%extremely efficient approach for solving high dimensional sparse
%regression problems. These
%algorithms work by cycling through all the predictors, at each step
%optimizing one parameter while holding all the other terms fixed.
Proposition~\ref{propL1} suggests Algorithm~\ref{alglin}, a
simple but very efficient coordinate descent algorithm for minimizing
(\ref{newopt2}) when $\rho_{\lambda_n}(t)=\lambda_nt$.

\begin{center}
\begin{algorithm}[b] \caption{\textbf{Linear FAR algorithm}}
\label{alglin}
%\noindent{\mathbf Linear FAR Algorithm}
%
\begin{enumerate}
\item[0.] Initialize~$\hat{\bfeta}_j={\mathbf0}$ and $S_j=\Theta
_j (\Theta_j^T\Theta_j )^{-1} \Theta_j^T$, for $j\in\{
1,\ldots,p\}$.
%\end{enumerate}
%For each~$j\in\{1,...,p\}$,
%\begin{enumerate}
%
\item[1.] Fix all $\hf_k$ for $k\ne j$. Compute the residual vector
$\bR_j = \bY-
\sum_{k\ne j} \hf_k$.
\item[2.] Let $\widehat\bP_j = S_j \bR_j$ represent the unshrunk
estimate for $\bff_j$.\vspace*{2pt}
%and $\phi_j = \rho'(\left\Vert\hf_j\right\Vert)$
% where $\hf_j$ represents the most recent estimate for $\bff_j$.
%
\item[3.] Let $\hf_j = \alpha_j \widehat\bP_j$ where $\alpha_j =
(1-\lambda_n\sqrt n/\llVert \widehat\bP_j\rrVert )_+$ is a
shrinkage parameter.\vspace*{3pt}
\item[4.] Center $\hf_j \leftarrow\hf_j - \operatorname{mean}(\hf_j)$.\vspace*{2pt}
\item[5.] Repeat steps~1 through 4 for $j=1,2,\ldots, p$ and iterate until
convergence.
\end{enumerate}
\end{algorithm}
\end{center}

We repeat this algorithm over a grid of values for $\lambda$, using
the previous values for the $\hat{\bfeta}_j$'s to initialize the
parameters for the new $\lambda$.
Since the parameters change very little for a small change in $\lambda
$, the algorithm generally converges very rapidly. Note that the
$S_j$'s only need to be computed
once for all values of $\lambda$ so the computation at each step of
the algorithm is extremely fast. In addition, it is clear from
Proposition~\ref{propL1} that
(\ref{newopt2}) will decrease at each step. This approach has the
advantage of
decomposing the estimation of $\hf_j$ into two simple, and separate,
steps. First, compute the unshrunk estimate $\widehat\bP_j$ and second,
apply the shrinkage factor $\alpha_j$. When $\alpha_j=0$ then the $j$th
predictor is absent from the model. Our FAR algorithm has similarities
to the SpAM algorithm \cite{ravikumar1} but SpAM cannot model
functional data.

For a general penalty function, $\rho_{\lambda_n}(t)$, we use the
local linear approximation method proposed in \cite{ZL08} to solve
(\ref{newopt2}). The penalty function can be approximated as $\rho
_{\lambda_n}(\llVert \bff\rrVert /\sqrt
n)\approx
\rho'_{\lambda_n}(\llVert \bff^*\rrVert /\sqrt n)\llVert \bff
\rrVert /\sqrt n+C$, where
$\bff^*$ is some vector that is close to $\bff$ and $C=\rho_{\lambda
_n}(\llVert \bff^*\rrVert /\sqrt
n)-\rho_{\lambda_n}'(\llVert \bff^*\rrVert /\sqrt n)\llVert \bff
^*\rrVert /\sqrt n$ is a
constant. Hence, the only required change to the FAR algorithm for
optimizing over general penalty
functions is to replace the calculation of $\alpha_j$ in step~3 by
\[
\alpha_j = \biggl(1-\rho'_{\lambda_n}\biggl(
\frac{1}{\sqrt{n}}\llVert\hf_j\rrVert\biggr)\sqrt n/\llVert\widehat
\bP_j\rrVert\biggr)_+,
\]
where $\hf_j$ represents the most recent estimate for $\bff_j$. The
initial estimate of $\hf_j$ can be obtained by using the Lasso
penalty. This simple
approximation allows the FAR algorithm to be easily applied to a wide
range of penalty functions.

%s2.3 #&#
\subsection{Nonlinear FAR}
\label{nlfar}

We now consider the more general nonlinear setting (\ref{simfar})
where $g_j(x)$ is estimated as part of the fitting process. Since
$\beta_j(t)$ corresponds to a direction that we project $X_{ij}(t)$
into we impose the constraint $\llVert \beta_j\rrVert _2=1$. Note
that $\beta_j$
are still not uniquely identifiable, however, our focus is on
estimating the regression functions, $f_j$, rather than the index functions.
%
%As in Section~\ref{linFARsec}, we assume that $\beta_j(t)$ and
%$X_{ij}(t)$
%are well approximated by an orthogonal $q$-dimensional basis $\bb(t)$
%such that $\beta_j(t) \approx\bb(t)^T\bfeta_j$ and
%$X_{ij}(t) \approx\bb(t)^T\btheta_{ij}$.
We assume that $g_j(x)$ can be well
approximated by a $d_n$-dimensional basis $\bh(x)$ such that $g_j(x)
\approx
\bh(x)^T \bxi_j$. Using this basis, representation (\ref{nonlinfar}) can
be expressed as
%
%e13 #&#
\begin{equation}
\label{nonlinfar2} \frac{1}{2n} \Biggl\llVert\bY-\sum
_{j=1}^p H_j\bxi_j\Biggr
\rrVert^2 + \sum_{j=1}^p
\rho_{\lambda_n} \biggl(\frac{1}{\sqrt n}\llVert H_j
\bxi_j\rrVert\biggr),
\end{equation}
where\vspace*{2pt} $H_j$ is an $n$ by $d_n$ matrix who's $i$th row is given by
$\bh(\btheta_{ij}^T \bfeta_j )^T$.

%\subsubsection{Nonlinear FAR Algorithm}

We use an iterative algorithm to approximately minimize (\ref
{nonlinfar2}) over $\bxi_j$ and $\bfeta_j$. First, given current
estimates for the $\bfeta_j$'s we minimize (\ref{nonlinfar2}) over
$\bxi_j$. Second, given current estimates for the $\bxi_j$'s we
minimize the sum of squares term
%
%e14 #&#
\begin{equation}
\label{etaopt} \sum_{i=1}^n
\Biggl(Y_i-\sum_{j=1}^p\bh
\bigl(\btheta_{ij}^T \bfeta_j
\bigr)^T\hbxi_j \Biggr)^2
\end{equation}
over $\bfeta_j$. Note that we do not include the penalty $\rho
_{\lambda_n}$ when estimating $\bfeta_j$ because the $\bfeta_j$'s
are providing a direction in which to project $X_{ij}(t)$ so are
constrained to be norm one. Hence, applying a shrinkage term would be
inappropriate.

Formally, the nonlinear FAR algorithm can be summarized as follows (Algorithm \ref{algnonlin}).
%
%\begin{center}
%
\begin{algorithm}[t] \caption{\textbf{Nonlinear FAR algorithm}}
\label{algnonlin}
%\noindent{\mathbf Linear FAR Algorithm}
%
\begin{enumerate}
\item[0.] Initialize\vspace*{2pt}~$\hat{\bfeta}_j$ for $j\in\{1,\ldots,p\}$
using the linear FAR algorithm.
%\end{enumerate}
%For each~$j\in\{1,...,p\}$,
%\begin{enumerate}
%
\item[1.] Compute $\widehat H_j$ using the current estimates for $\bfeta_j$.
\item[2.] Estimate $\bxi_j$ for $j\in\{1,\ldots,p\}$ by minimizing
(\ref{nonlinfar2}) given the current values of $\widehat H_j$.
\item[3.] Conditional on the $\hbxi_j$'s from step~2, estimate the
$\bfeta_j$'s by minimizing (\ref{etaopt}).

\item[4.] Repeat steps 1~through 3 and iterate until convergence.
\end{enumerate}
\end{algorithm}
%
%\end{center}

One of the appealing aspects of this approach is that, for fixed $\widehat
H_j$, (\ref{newopt2}) and (\ref{nonlinfar2}) are equivalent so
estimation of the $\bxi_j$'s in step 2 can be achieved using the
linear FAR algorithm from Section~\ref{seclinearfar}. Minimization of
(\ref{etaopt}) in step~3 can be approximately achieved using a
first-order Taylor series approximation of $g_j(x)$. We provide details
on this minimization and on computing initial values for the $\bfeta
_j$'s in Appendix~\ref{stepA}.

Potentially one could compute the nonlinear FAR algorithm for each
possible~$\lambda$. However, we have found that a more efficient
approach is to compute initial estimates for $\bfeta_j$, minimize
(\ref{nonlinfar2}) over $\bxi_j$ for each possible value of $\lambda
$, choose the $\bxi_j$'s corresponding to the ``best'' value of
$\lambda$, estimate the $\bfeta_j$'s for only this one set of
parameters, and iterate. This approach means that, for each iteration,
the minimization of (\ref{etaopt}) only needs to be performed for a
single value of $\lambda$. The choice of $\lambda$ can be made using
a variety of methods, as discussed in the next section.

%
%{\bolds
%\noindent{\mathbf Algorithm}
%\begin{enumerate}
%\item[A.] Given an estimate $\hbfeta$ for $\bfeta= (\bfeta_1\t,
%\cdots, \bfeta_p\t)\t$, minimize
% $l_\lambda(\bxi\left\vert\hat{\bfeta})$ (\ref{nonlinfar2}) over $
%\bxi_1,
%\ldots, \bxi_p$ using the linear FAR algorithm to get estimate $\hbxi
%= (\hbxi_1\t,\cdots, \hbxi_p\t)\t$.
%
%\item[B.] Conditional on $\hbxi= (\hbxi_1\t,\cdots, \hbxi_p\t)\t$
%from Step A, compute the $\bfeta_j$'s by minimizing
%\begin{equation}
%\label{etaopt}
%Q(\bfeta\right\vert\hbxi)=\sum_{i=1}^n \Big(Y_i-\sum_{j=1}^p\bh\left(
%\btheta_{ij}^T
% \bfeta_j\right)^T\hbxi_j\Big)^2.
%\end{equation}
%Our approach for minimizing (\ref{etaopt}) is provided in Appendix~
%\ref{stepA}.
%
%\item[C.] Iterate between Step A. and Step B. until the value of $
%\hbfeta$ converges or the maximum number of iterations is reached.
%\end{enumerate}

%In practice, Step A. is implemented over a grid of tuning
%parameters, $\lambda_1,\ldots, \lambda_T$, which makes the algorithm
%slow to converge. To speed it up, starting from the second iteration,
%we can fix the value of $\lambda$ at the one chosen in Step A. of the
%first iteration. This simplification works well in our simulations and
%real data analysis. We discuss our approach for selecting initial
%values $\hbfeta_j$'s in Appendix~\ref{initialvalues}.

%s2.4 #&#
\subsection{Selecting tuning parameters}

Both the linear and nonlinear versions of FAR require choosing the tuning
parameter, $\lambda$. As with all penalized regression methods, there are
several possible methods one could adopt. Popular approaches include,
BIC, AIC or cross-validation. The BIC and AIC methods require the
calculation of the effective degrees of freedom. For the Lasso, it has been
shown that an unbiased estimate for this quantity is the number of
nonzero coefficients \cite{zou2}. One could potentially use the
same value for
FAR. However, given FAR's more complicated structure it is not clear that
this is still an appropriate estimate. Computing the effective degrees of
freedom for FAR is a topic for future research. For our simulations and
one real data example, we
selected $\lambda$
using a separate validation data set. For the other real data example,
we selected $\lambda$ using the 20-fold cross-validation method, since
there were not enough data points to be used as validation data.
%cross-validation method. %for our real data examples and a separate
%validation data set for our simulations. This approach seemed to work
%well on the problems that we
%examined. {\mathbf How did we do it for simulation data?}

%\subsection{Accelerating the FAR Algorithm}

%s3 #&#
\section{Theory}
\label{theorysec}

%s3.1 #&#
\subsection{Linear theory}

Denote by $\Mfrak_0 = \{j: \beta_j(t)\neq0, 1\leq j\leq p\}$ the set
of true functional predictors and let $s_n$ represent the cardinality
of $\Mfrak_0$.
By minimizing the FAR criterion (\ref{newopt2}), we aim to identify
the set $\Mfrak_0$ and accurately estimate functions $\beta_j(t)$ for
$j\in\Mfrak_0$. In this section,
we discuss the theoretical properties of FAR in the setting where the
$f_j$'s are linear functions, that is, $f_j(X_{ij}) = \int
_0^1X_{ij}(t)\beta_j(t)\,dt$. In particular,
we present two theorems, both of which are conditional on the observed
predictors, $X_{ij}(t)$, $i=1,\ldots, n$, $j = 1, \ldots, p$.
Theorem~\ref{T2} concerns FAR's model selection properties. We show
that, with probability tending to one, FAR can remove all noise predictors
from the fitted model. Theorem~\ref{T2} also places an error bound on
the estimated $f_j$'s under the vector infinity norm, where $j\in
\Mfrak_0$. Our second result, Theorem~\ref{T3}
shows the asymptotic normality of the estimator.
%Proofs of all the results in this section are provided in the appendix.

In order to prove these results, we make two sets of assumptions. The
first set of conditions relates to the level of accuracy in our basis
approximations of $X_{ij}(t)$ and
$\beta_j(t)$. The second set of conditions concerns the shape of the
penalty function, the strength of the signal and the correlation
structure of the predictors. Explicit conditions %and proofs of all the
%results in this section,
can be found in Appendix~\ref{linearCon}.

Let $\bfeta_0 =
(\bfeta_{0,1}, \ldots, \bfeta_{0,p}) \in R^{pq_n}$ with $\bfeta
_{0j}$ representing the true coefficient vector in the basis
representation $f_j(X_{ij}) = \btheta_{ij}^T\bfeta_{0j} + e_{ij}$.
%In view of (\ref{e005}), we can estimate $\Mfrak_0$ by solving the
%above regularization %problem (\ref{e005}).
For\vspace*{1pt} any index set $S \subset\{1,\ldots, p\}$, we use $\bfeta_S$ to
denote the vector formed by stacking vectors $\bfeta_j$, $j\in S$ one
underneath each other, and
$\Theta_S$ to denote the matrix formed by stacking the matrices
$\Theta_j$, $j\in S$ one after another. Moreover, we standardize each
column of $\Theta$ such that they all have $\ell_2$-norm $\sqrt{n}$.
Theorem~\ref{T2} below shows that FAR possesses the oracle property
for model selection.

%th1 #&#
\begin{theorem}\label{T2} Assume that $q_n + \log p = O(n\lambda_n^2)
$, $\lambda_n n^{\alpha}q_n\sqrt{s_n} \rightarrow0$, and $\log
(pq_n) = o(n^{1-2\alpha}s_n^{-1}q_n^{-2})$ with $\alpha$ defined in
Condition~\ref{as2}\textup{(B)}. Further assume that $s_nq_n^{-2} = o(\lambda
_n)$, then under Conditions~\ref{as3} and~\ref{as2}, with probability
tending to 1 as $ n
\rightarrow\infty$, there exists a local minimizer $\hbfeta$ of
(\ref{newopt2}) such that:
\begin{longlist}[(2)]
\item[(1)] $\hbfeta_{\Mfrak_0^c}=0$,
\item[(2)] $\llVert \hbfeta_{\Mfrak_0} - \bfeta_{0\Mfrak_0}\rrVert
_\infty\leq
c_0^{1/2}n^{-\alpha}q_n^{-1/2}$,
\end{longlist}
where $\llVert \cdot\rrVert _\infty$ stands for the infinity norm of
a vector.
\end{theorem}

Although Theorem~\ref{T2} is on a local minimizer of the linear FAR
criterion (\ref{newopt2}), it has been proved by \cite{loh13} that
any local minimizer will fall within statistical precision of the true
parameter vector under appropriate conditions on the penalty function.
Part 2 of Theorem~\ref{T2} concerns the approximation accuracy of the
basis coefficients rather than the functions themselves. However, the
result extends naturally.
Denote by $\hat\bff_j = \Theta_j\hbfeta_j$ and $\bff_{0j} =
(f_j(X_{j1}), \ldots, f_j(X_{jn}))^T$, respectively, the estimated and
true values of the $j$th functional
component, both evaluated at the $n$ training data points. Then the
corollary below follows immediately from Theorem~\ref{T2} and
Condition~\ref{as3}.

%co1 #&#
\begin{corollary}\label{cor1} Suppose the conditions in Theorem~\ref
{T2} are satisfied. Then with probability tending to 1 as $n\rightarrow
\infty$, there exists a FAR
estimate such that $\hbff_j = 0$ for $j\notin\Mfrak_0$, and
\[
\max_{j\in\Mfrak_0}\frac{1}{\sqrt{n}}\llVert\hbff_j -
\bff_{0j}\rrVert_2 \leq C_2n^{-\alpha},
\]
%
%\[
%\max_{j\in\Mfrak_0}\left\Vert\hat f_j -f_{j}\right\Vert_{L_2(\widehat
%F_j)}= \max_{j\in
%\Mfrak_0}\frac{1}{\sqrt{n}}\left\Vert\hbff_j - \bff_{0j}\right\Vert
%\leq C_2n^{-\alpha}
%\sqrt{\log n},
%\]
where $C_2$ is some positive constant.
%, $\hat f_j(x) = \int_0^1x(t)\widehat\beta_j(t)\,dt$ with $\widehat\beta_j(t) =
%\bb^T(t)\hbfeta_j$, $x(t)$ is a given value of the $j$-th functional
%predictor, and $\widehat
%F_j$ is the corresponding empirical distribution of $x(t)$.
\end{corollary}

Theorem~\ref{T3} shows the asymptotic normality of the FAR estimators
that correspond to signal variables. As with Theorem~\ref{T2}, we
first provide the result for the
$\hbfeta_j$'s and then extend to the functions.

%th2 #&#
\begin{theorem} \label{T3}
Assume that the conditions in Theorem~\ref{T2} hold and in addition,
$\rho_{\lambda_n}'(a_n/2)=o(a_nn^{\alpha-\sfrac{1}{2}}s_n^{-1/2})$,
$\sup_{t\geq\sfrac{a_n}{2}}\rho_{\lambda_n}''(t)=O(n^{-1/2})$,
$s_{n}=o(n^{2\alpha})$ and $s_nq_n^{-2} = o(n^{-1/2})$. Then with
probability tending to 1
as $n \rightarrow\infty$, there exists a strict local minimizer
$\hbfeta$
of (\ref{newopt2}) such that $\hbfeta_{\Mfrak_0^c} = 0$ and
\[
\bc^T \bigl[\bigl(\Theta_{\Mfrak_0}^T
\Theta_{\Mfrak_0}\bigr)^{1/2}(\hbfeta_{\Mfrak_0} -
\bfeta_{0,\Mfrak_0}) + n\bigl(\Theta_{\Mfrak_0}^T\Theta
_{\Mfrak_0}\bigr)^{-1/2}\bv_{0,\Mfrak_0}\bigr] \toD N\bigl(0,
\sigma^2\bigr), %
\]
where $\bc\in\mathbf{R}^{q_ns_n}$ satisfies $\bc^T\bc= 1$ and
$\bv_{0,\Mfrak_0}$ is a vector formed by stacking the
vectors
$\bv_{0,k}=\rho'_{\lambda_n}(\frac{1}{\sqrt{n}}\llVert \Theta_k \bfeta
_{0,k}\rrVert )\frac{1}{\sqrt{n}} \frac{\Theta_k^T\Theta_k\bfeta
_{0,k}}{\llVert \Theta_k\bfeta_{0,k}\rrVert }$,
$k\in\Mfrak_0$ underneath each other.
\end{theorem}

Let $f_{0j}^* = {\btheta_j^*}^T\bfeta_{0j}$ and $\hat f_{j}^* =
{\btheta^*_j}^T\hbfeta_{j}$, with $\btheta_j^*\in R^{q_n}$ the
coefficient vector when projecting a given
new observation, $X^*_j(t)$, onto the basis function, $\bb(t)$. Then
as $q_n$ increases, $f^*_{0j}$ better approximates $f_j(X_j^*)$ for
each fixed $j = 1,\ldots, p$.
Define $\bff_0^* = (f_{01}^*, \ldots, f_{0p}^*)^T$ and $\hf^* =
(\hat f_1^*, \ldots, \hat f_p^*)^T$. Taking $\bc=
(\Theta_{\Mfrak_0}^T\Theta_{\Mfrak_0})^{-1/2}\Theta^*\tilde\bc
_0$ with $\Theta^* = \operatorname{diag} ({\btheta_1^*}, \ldots,
{\btheta_{s_n}^*} )\in R^{(q_ns_n) \times s_n}$
in Theorem~\ref{T3} and $\tilde\bc_0$ a vector in $R^{s_n}$, we have
the following asymptotic normality of $\bff_0^*$.

%co2 #&#
\begin{corollary} Assume that the conditions in Theorem~\ref{T3} hold.
Then with probability tending to 1
as $n \rightarrow\infty$, there exists a FAR estimate such that
$\hbff_{\Mfrak_0^c}^* = 0$. Moreover,
\[
\tilde\bc^T_0\bigl[\hf_{\Mfrak_0}^* -
\bff_{0,\Mfrak_0}^* + n\bigl(\Theta^*\bigr)^T\bigl(
\Theta_{\Mfrak_0}^T\Theta_{\Mfrak_0}\bigr)^{-1}
\bv_{0,\Mfrak
_0}\bigr]\toD N\bigl(0, \sigma^2\bigr),
\]
where $\tilde\bc_0$ is a vector in $R^{s_n}$ satisfying $\tilde\bc
_0^T(\Theta^*)^T(\Theta_{\Mfrak_0}^T\Theta_{\Mfrak_0})^{-1}\Theta
^*\tilde\bc_0 = 1$, and $\bv_{0,\Mfrak_0}$
is defined in Theorem~\ref{T3}.
\end{corollary}

%s3.2 #&#
\subsection{Nonlinear theory}
\label{nonlintheorysec}

Throughout this section, we focus on the minimizer of the nonlinear
FAR criterion with the $\ell_1$ penalty function. We treat all the
predictors as deterministic. For identifiability purposes, we assume
that the true regression functions, $f_{0j}$, as well as the response
vector, are centered, that is, $\sum_{i=1}^n f_{0j}(X_{ij})=0$ and
$\sum_{i=1}^n Y_i=0$. As a result, the corresponding estimates,
$\hat f_j$, are automatically centered as well.
%and allow an intercept term in both the nonlinear far model, (
%\ref{simfar}), and the corresponding optimization criteria, (
%\ref{nonlinfar}) and~(\ref{nonlinfar2}).

We use cubic B-splines to approximate the true ``link'' functions, $g_{0j}$.
%For simplicity of the exposition, we will set the dimension of each
%such basis equal to the dimension of the orthonormal basis used in
%representing the predictor functions~$X_{ij}(t)$ and the index
%functions~$\beta_j(t)$. It follows from the proofs that such
%simplification does not hurt the results of this section under the
%assumptions that we impose. Thus, the dimensions~$q_n$ and~$d_n$,
%introduced earlier, will equal each other, and will be referred to
%as~$d_n$.
Given a candidate index vector $\bfeta_j$, the B-spline basis for
representing a candidate link function for the $j$th predictor is
constructed using uniformly placed knots on the interval $[\min_i
\bfeta^T_j\btheta_{ij},\max_i \bfeta^T_j\btheta_{ij}]$. The
corresponding row vector valued basis function is denoted by~$\bh
_{\bfeta_j,j}$. We denote by $\F_j^0$ the class of candidate
regression functions for the $j$th predictor. More specifically, $\F
_j^0=\{f(\cdot)=\bh_{\bfeta_j,j}(\bfeta_j^T\cdot)\bxi, \sum_{i=1}^n
f(\btheta_{ij})=0, \bfeta_j\in\RR^{q_n}, \bxi\in\RR
^{d_n}, \llVert \bfeta_j\rrVert =1\}$. If~$\tilde f$ and~$\check{f}$ belong
to~$\F_j^0$, we denote by $\llVert \tilde f-\check{f}\rrVert _n$ the $\ell_2$
distance between these two functions with respect to the empirical
probability measure corresponding to $\btheta_{1j},\ldots,\btheta
_{nj}$. More specifically, $\llVert \tilde f-\check{f}\rrVert
^2_n=n^{-1}\sum_{i=1}^n (\tilde f(\btheta_{ij})- \check{f}(\btheta
_{ij}))^2$. We\vspace*{1pt}
refer to the estimated regression function for the $j$th predictor as
$\hat f_j(\cdot)=\hat g_j(\hat{\bfeta}_j^T\cdot)$. The\vspace*{1pt}
corresponding true regression functions are referred to as $f_{0j}$. We
slightly abuse the notation and write $\llVert \hat f_j-f_{0j}\rrVert _n$ for
the $\ell_2$ distance between $\hat f_j$ and $f_{0j}$ with respect to
the empirical probability measure corresponding to the $j$th predictor:
\begin{eqnarray*}
\llVert\hat f_j-f_{0j}\rrVert^2_n&=&
\frac{1}n\sum_{i=1}^n \bigl[
\hat f_j(\btheta_{ij})- f_{0j}(X_{ij})
\bigr]^2
\\
&=&\frac{1}n\sum_{i=1}^n \biggl[
\hat g_j \bigl(\hat{\bfeta}_j^T
\btheta_{ij} \bigr)- g_{0j} \biggl(\int_0^1
\beta_j(t)X_{ij}(t)\,dt \biggr) \biggr]^2.
\end{eqnarray*}

As before, we write~$\Mfrak_0$ for the index set of the signal
predictors, that is, $\Mfrak_0 =\{j: 1\le j\le{p_n}, f_{0j}\ne0\}$.
Note that this set depends on~$n$, but we will refrain from using an
additional subscript for simplicity of the notation. We use~$\widehat
\Mfrak_n$ to denote the corresponding estimated set, $\{j: 1\le j\le
{p_n}, \hat f_{j}\ne0\}$. Let~$s_n=\llvert \Mfrak_0\rrvert $. A
\textit{universal constant} is interpreted as a constant that does not
depend on~$n$ or any of the other parameters that appear in the
corresponding expression. Given expressions~$E_1$ and~$E_2$, we
use~$E_1\gtrsim E_2$ to mean that there exists a positive universal
constant~$c$, such that $E_1\ge cE_2$. We write $E_1\asymp E_2$ when
both $E_1\gtrsim E_2$ and $E_2\gtrsim E_1$ are satisfied.

The results provided below establish the rate of convergence for the
estimated regression functions. To derive these results, we impose a
number of regularity conditions on the components of the FAR model. We
also impose a version of the \emph{compatibility condition}, which is
commonly used in high-dimensional additive models \cite
{meier1,Bbuhlmann1}. The proofs, as well as a more detailed
discussion of the conditions, are provided in Appendix~\ref{appnonlimthm}.

%th3 #&#
\begin{theorem} \label{nonlinthm}
Suppose that Conditions~\ref{cond3} and~\ref{as4} are satisfied. Let
$q_n\gtrsim d_n\gtrsim\log\log n$. Then there exists a universal
constant $c$, such that for $\lambda_n\ge
c(n^{-1/2}q_n^{1/2}+n^{-1/2}\sqrt{\log{p_n}})$, the following bound
holds with probability tending to one, as $n$ tends to infinity:
%
%e15 #&#
\begin{equation}
\label{nonlinthmbnd} \sum_{j=1}^{{p_n}} \llVert
\hat f_j-f_{0j}\rrVert_n= O
\bigl(s_n\lambda_n+s_nd_n^{-2}+s_n^2n^{1/2}d_n^{-4}q_n^{-1/2}
\bigr).
\end{equation}
\end{theorem}

The following corollary focuses on the choice of~$q_n$ and $d_n$ that
yields the fastest rate of convergence. Note that the case
$q_n/d_n=o(1)$ is not covered in the statement of Theorem~\ref
{nonlinthm}. However, it follows from the proof of the theorem that
such settings correspond to an error bound that is inferior to the one
presented below.

%co3 #&#
\begin{corollary}
\label{nonlinratecor}
Suppose that Conditions~\ref{cond3} and~\ref{as4} are satisfied. Let
$q_n\asymp d_n\asymp(s_nn)^{1/5}$. Then\vspace*{1pt} there exists a universal
constant $c$, such that for $\lambda_n \ge
c(s_n^{1/10}n^{-2/5}+n^{-1/2}\sqrt{\log{p_n}})$, the following bound
holds with probability tending to one, as $n$ tends to infinity:
\[
\sum_{j=1}^{{p_n}} \llVert\hat f_j-f_{0j}\rrVert_n= O(s_n
\lambda_n).
\]
\end{corollary}

We now turn to the variable selection properties of the nonlinear FAR
estimator. Methods that use~$\ell_2$ regularization are known to typically
produce models containing a large number of noise predictors (\cite
{Bbuhlmann1}, Chapter~7, e.g.). To alleviate this problem, we
follow the popular approach of thresholding the initial estimator. We
define the thresholded FAR estimator as follows: $\tilde f_j=\hat f_j I{\{\llVert \hat f_j\rrVert _n> \lambda_n\}}$,
$j=1,\ldots
,p_n$. Note that the threshold parameter is taken equal to the tuning
parameter~$\lambda_n$, which is used to compute the initial
estimators, $\hat f_j$. Thus, we do not introduce any new tuning
parameters at the thresholding stage. Let~$\widetilde\Mfrak_n$ denote
the index set of the corresponding nonzero regression function
estimates, $\{j: 1\le j\le{p_n}, \tilde f_{j}\ne0\}$. Recall that
$s_n=\llvert \Mfrak_0\rrvert $. The next result provides bounds for
the estimation
error of the thresholded FAR approach and for the corresponding number
of selected predictors.

%th4 #&#
\begin{theorem} \label{threshthm}
Under\vspace*{1pt} the assumptions of Corollary~\ref{nonlinratecor}, there exists
a universal constant $c$, such that for $\lambda_n \ge
c(s_n^{1/10}n^{-2/5}+n^{-1/2}\sqrt{\log{p_n}})$, the following bounds
hold with probability tending to one, as $n$ tends to infinity:
\begin{eqnarray*}
\llvert\widetilde\Mfrak_n\rrvert &=& O(s_n)\quad
\mbox{and}
\\
\sum_{j=1}^{{p_n}} \llVert\tilde f_j-f_{0j}\rrVert_n &=& O(s_n
\lambda_n).
\end{eqnarray*}
\end{theorem}

Now consider the case where the components of the FAR model do not
depend on~$n$. More specifically, suppose that the number of signal
predictors, $\llvert \Mfrak_0\rrvert $, and the signal regression
functions, $\{
f_{0k}\}_{k\in\Mfrak_0}$, are fixed and do not change with~$n$. The
estimation error bound in Theorem~\ref{threshthm} implies that, with
probability tending to one, our estimator has zero false negatives,
while the number of false positives stays bounded. This variable
selection result can be strengthened by increasing the threshold
from~$\lambda_n$ to~$\tau\lambda_n$, for a sufficiently large~$\tau
$. The next corollary demonstrates that the corresponding thresholded
estimator can correctly recover the index set of the relevant predictors.

%co4 #&#
\begin{corollary}
\label{threshcor}
Suppose that the components of the FAR model do not depend on~$n$.
Suppose also that the assumptions of Corollary~\ref{nonlinratecor}
are satisfied. Then there exist universal constants~$\tau_0$ and~$c$,
such that, provided $\tau\ge\tau_0$, $\lambda_n \ge
c(s_n^{1/10}n^{-2/5}+n^{-1/2}\sqrt{\log{p_n}})$ and $\lambda
_n=o(1)$, we have
\[
\widetilde\Mfrak_n=\Mfrak_0,
\]
with probability tending to one, as~$n$ goes to infinity.
\end{corollary}

%s4 #&#
\section{Simulations}
\label{simsec}

In this section, we compare the performance of FAR to several alternative
linear and nonlinear functional approaches in a series of simulation
studies. We consider the linear setting in Section~\ref{linear}, while
Section~\ref{nonlinear} contains our nonlinear results.

%s4.1 #&#
\subsection{Linear additive models}
\label{linear}

We\vspace*{1pt} first generated the functional predictors, $X_{ij}(t)$, from a
4-dimensional Fourier basis $\bb(t)= (1, \sqrt{2}\sin(\pi t), \sqrt
{2}\sin(2\pi t)$, $ \sqrt{2}\sin(3\pi t))^T$, plus an error term:
\begin{eqnarray*}
&& X_{ij}(t_k) = \bb(t_k)^T
\btheta_{ij}+w_{ijk}, \qquad w_{ijk} \sim N\bigl(0,
\sigma_x^2\bigr), \qquad\btheta_{ij} \sim N(0,
I),
\end{eqnarray*}
where $\sigma_x = 0.5$, and each predictor was observed at $200$ equally
spaced time points, $0=t_1, t_2,\ldots, t_{200}=1$.
%Each data set contained $n=60$ observations.
The basis
coefficients, $\btheta_{ij}$, and the error terms, $w_{ijk}$, were all
sampled independently from each other.
The first $s_n$ coefficient functions, $\beta_1(t), \ldots, \beta
_{s_n}(t)$, were
also generated, from the same basis function, $\beta_j(t)=\bb
(t)^T\bfeta_j$,
while the remaining $p-s_n$ predictors were noise variables with
$\beta_j(t)=0$. For each $j=1,\ldots, s_n$, the coefficient vector
$\bfeta_j$ were first independently generated from a multivariate
standard normal distribution and then rescaled to have $\ell_2$ norm
equal to 1. The responses were then generated from (\ref{far}) with
$f_j(x)$ computed using (\ref{linearf}). We tested a total of six
linear settings corresponding to different numbers of observations,
predictors and noise levels.

%In this subsection we generated data from the linear FAR model with $
%\beta_0 = 0$
%\begin{eqnarray*}
%Y_i = \sum_{j=1}^p\int_0^1X_{ij}(t)\beta_j(t)\,dt + \varepsilon_i, i=1,
%\cdots, n,
%\end{eqnarray*}
%$\varepsilon_i \sim_{i.i.d.} N(0, \sigma_y^2)$.

To ensure a fair real world comparison, where the true functional form
of $\beta_j(t)$ would be unknown, we implemented the linear version of
FAR using an orthogonal cubic spline basis, rather than the true
Fourier basis. We tested FAR using both the SCAD \cite{fan2} and the
Lasso penalty functions but found that the former penalty generally
gave superior predictive ability so only report the SCAD results here.
%We also generated unshrunk estimates (FAR*) by first selecting a
%subset of predictors using
%the regular FAR methodology, and then fitting this subset of variables
%using an unpenalized version of FAR. This approach can reduce the
%overshrinkage often exhibited by the Lasso penalty and is analogous to
%the LARS/OLS hybrid discussed in \cite{efron3}.
%
We compared FAR to three competing methods. The first was a functional principal
components analysis (FPCA) based approach produced by decomposing the
predictors into functional principal components, selecting the first
$K$ components
and finally using the resulting PCA scores to fit linear regression
models to the response. Since only $s_n$ of the predictor functions
were associated with the response, we fit the linear regressions to the
FPCA scores using the group SCAD penalty function to produce sparse
fits, where the $K$ principal components for each predictor were
grouped together.
%The FPCA estimates were produced using the standard smoothing
%approach; compute the basis matrix $B$ associated with
%the spline basis $\bb(t)$, individually estimate the
%$\btheta_{ij}$'s using the least squares fit, $\widehat\btheta_{ij} =
%(B^TB)^{-1} B^T \bX_{ij}$, and finally apply standard PCA to the
%resulting $\widehat\btheta_{ij}$'s.
%%As with FAR we examined both the regular shrunk coefficient estimates
%(FPCA) and the unshrunk estimates (FPCA*).

Our second approach involved implementing the additive modeling method
(ADD) of \cite{ferraty3}. ADD fits an additive model with the same
general form as (\ref{far}). A key difference relative to FAR is that
ADD uses a kernel based fitting method and a forward selection
procedure to iteratively add functional predictors to the model. The
final method, SIR, is described in \cite{amato1}. This method first
computes the wavelet coefficients on a single predictor function, then
applies the SIR \cite{li1} dimension reduction method to the
resulting coefficients, and finally a linear regression is fit using
the reduced dimensions as the predictors. This approach is not designed
for multiple predictor functions so we adapted it by computing the
reduced dimensions marginally for each predictor and then performing a
multiple linear regression on all the resulting dimensions.

%Let $\gamma$ represent the proportion of variation in $
%\int_0^1X_{ij}(t)\beta_j(t)\,dt$ that is explained by the first principal
%component of $X_{ij}(t)$. Values of $\gamma$ close to one correspond
%to very favorable situations for the FPCA methods since almost all the
%information about the response is contained in the first
%principal component of $X(t)$. We tested a total of six linear
%settings corresponding to $\gamma\approx69\%$, $90\%$ or $99\%$, $
%\sigma=1$ or $2$, and, $p=10, 50$ or $200$. The $\gamma=90\%$ and $
%\gamma=99\%$ settings respectively corresponded to favorable and
%extremely favorable situations for FPCA while $\gamma=69\%$
%represented a more balanced scenario.

The tuning parameters for the various methods were chosen by minimizing
prediction error on a separately generated validation data set with
identical characteristics to the training data. FAR had two tuning
parameters; $\lambda$ and the dimension of the orthogonal cubic spline
basis for fitting $\beta_j(t)$. We fitted FAR separately for each
possible basis dimension, and then selected the value (between $5$ and
$10$) which gave the smallest prediction error on the validation set.
The FPCA method had two tuning parameters; $\lambda$, the penalty
level for the group SCAD fit, and $K$, the number of principal
components used for each predictor. We used the same value of $K$ for
all predictors. To select $K$, we first identified a number $K_{\max}$
such that the first $K_{\max}$ scores of each predictor express at
least 99\% of the total variation of this predictor, and then selected
$K$ as the value (between 1 and $K_{\max}$) which minimized prediction
error on the validation data. The SIR method had one tuning parameter;
the number of directions into which each predictor was projected. We
considered up to $4$ directions for each predictor, and selected the
number of directions as the one with the lowest prediction error on the
validation set.

For each simulation setting, we fitted each method to $100$ different
training sets and recorded the false positive rate (FPR), false
negative rate (FNR), average prediction error on a separate test data
set (Mean PE) and the standard error of the mean PE (SE PE). The FPR
records the fraction of noise predictors incorrectly included in the
model while the FNR corresponds to the fraction of signal variables
incorrectly excluded. The simulation results are summarized in
Table~\ref{tablin}. Prediction errors that were either the best or
were not statistically worse than the best result are shown in bold
font. Note that because of the extremely computationally intensive
nature of the ADD and SIR methods it was not feasible to compute fits
for $p$ larger than about $10$. In fact, in the $p=600$ and 2000
settings the FPCA, ADD and SIR comparison methods were all too slow to
implement, and thus we only report the results for FAR. In terms of
prediction error, FAR was superior to all of the competing methods in
most simulation settings. The FPCA method was the best competitor
followed by SIR and finally ADD.
The only setting where FPCA was superior was the situation where
$\sigma_y=2$ and $p=100$, which had high noise and high
dimensionality. For the ultra-high dimensional setting of $p={}$2000,
FAR still does a reasonably good job in variable selection. Note that
when fitting FAR, since each functional predictor is approximated using
a spline basis, the dimensionality in the linear FAR criterion is in
fact much higher than $p$. For example, if a 5-dimensional spline basis
is used, the dimensionality is in fact $5p ={}$10,000.

%t1 #&#
\begin{table}%[t]
\tabcolsep=0pt
\caption{Comparison of FAR to three alternative methods in five linear
simulation settings}\label{tablin}
\begin{tabular*}{\tablewidth}{@{\extracolsep{\fill}}@{}lccccc@{}}
\hline
& & \multicolumn{1}{c}{\textbf{FAR}} & \multicolumn{1}{c}{\textbf{FPCA}} & \multicolumn{1}{c}{\textbf{ADD}} & \multicolumn{1}{c@{}}{\textbf{SIR}}\\
\hline
$n=60$ & FN & 0.0000 & 0.0000 & NA & NA \\
$p=10$ & FP & 0.0250 & 0.1067 & NA & NA \\
$s_n=4$ & Mean PE & \textbf{1.4834} & 1.6558 & 2.6474 & 2.3318 \\
$\sigma_y=1$ & SE PE & 0.0285 & 0.0274 & 0.0298 & 0.0275
\\[6pt]
$n=60$ & FN & 0.0225 & 0.005\phantom{0} & NA & NA \\
$p=10$ & FP & 0.05\phantom{00} & 0.1633 & NA & NA \\
$s_n=4$ & Mean PE & \textbf{2.6805} & 2.7979 & 3.3462 & 6.2968 \\
$\sigma_y=2$ & SE PE & 0.0267 & 0.0264 & 0.0296 & 0.0857
\\[6pt]
$n=80$ & FN & 0.0067 & 0.1917 & & \\
$p=100$ & FP & 0.0743 & 0.0454 & & \\
$s_n=6$ & Mean PE & \textbf{2.0176} & 3.5502 & & \\
$\sigma_y=1$ & SE PE & 0.0548 & 0.0353 & &
\\[6pt]
$n=80$ & FN & 0.0483 & 0.0067 & & \\
$p=100$ & FP & 0.1896 & 0.1569 & & \\
$s_n=6$ & Mean PE & 3.7051 & \textbf{3.3250} & \\
$\sigma_y=2$ & SE PE & 0.0548 & 0.0353 &
\\[6pt]
$n=100$ & FN & 0.0700 & &&\\
$p=600$ & FP & 0.0432&&& \\
$s_n=8$ & Mean PE & 3.6423 & &&\\
$\sigma_y=1$ & SE PE & 0.0910 &&&
\\[6pt]
$n=100$ & FN & 0.1925 & &&\\
$p=2000$ & FP & 0.0171&&& \\
$s_n=8$ & Mean PE & 4.6422 & &&\\
$\sigma_y=1$ & SE PE & 0.0871 &&& \\
\hline
\end{tabular*}
\end{table}

%s4.2 #&#
\subsection{Nonlinear models}
\label{nonlinear}
We examined three different simulation settings with the responses
generated from the nonlinear model (\ref{simfar}). The standard
deviation, $\sigma_x$, the predictors, $X_{ij}(t)$, and coefficient
curves, $\beta_j(t)$, were all produced in an
identical fashion to the linear setting. To produce a sparse
relationship between the predictors and the response, we set $g_j(x)=0$
for $j=3,4,\ldots, p$. The remaining two curves were chosen as
$g_1(x_1)=x_1$ and $g_2(x_2)=-x_2+\sin(x_2)$. Note that these
functions were not generated from a B-spline basis so the FAR fit
contains bias in the estimates for both $\beta_j(t)$ and $g_j(x)$; a
real world situation where the data is unlikely to exactly correspond
to the FAR model. The sample size was fixed at $n=100$, and the model
errors were independently generated from a Gaussian distribution with
mean zero and standard deviation $\sigma=0.5$.

We compared the nonlinear version of FAR to the same three competing
methods as in the linear setting. However, to account for the
nonlinear relationships between the response and predictors, we
implemented FPCA by applying the SpAM method \cite{ravikumar1} to
the principal component scores. SpAM essentially fits a penalized
version of Generalized Additive Models (GAM), allowing for automatic
variable selection in a nonlinear but additive regression situation.
We adapted SpAM slightly to implement a group penalization where all
$K$ PCs for a given predictor were penalized together. The SIR method
was still implemented using the linear regression approach from the
previous section while the kernel approach of ADD already produced a
nonlinear fit so these last two methods did not require any
adaptations to the new setting. In each simulation, we again fit the
methods to $100$ separate data sets and used a separate validation data
set, with identical
characteristics to the training data, to select the tuning parameters.
The nonlinear setting increased by one the number of tuning parameters
for the FAR and FPCA methods; $d$, the basis dimension for $g_j(x)$.
For both methods, we chose $d$ by computing the validation error rates
for values between $5$ and $10$, selecting the optimal value and then
using this dimension to compute $g_j(x)$. To reduce the computational
cost for FAR, we selected $q$, the dimension of the spline basis for
$\beta_j(t)$, as the value (between $5$ and $10$) which gave the best
hold out accuracy on the predictors in the validation set. In
particular, we held out $20\%$ of each predictor's time points,
computed the least squares fit to the remaining time points for each
possible basis dimension, and then selected the value of $q$ which gave
the lowest error rate on the held-out points.

%\begin{table}[t]
%\begin{center}
%{\small\begin{tabular}{r\right\vert r\right\vert rrrr}
%\hline\hline
% & & FAR & FPCA & ADD & SIR \\
% \hline
%$n=100$ & FN & 0.0000& 0.0000 &&\\
%$p=5$ & FP & 0.1033& 0.1300 &&\\
%$\sigma_y=0.5$ & Mean PE & {0.9762}& 1.3108& 1.7408& \textbf{0.8688}\\
% & SE PE & 0.0129& 0.0174& 0.0074& 0.0049\\
%\hline
%$n=100$ & FN & 0.0000 & 0.0000 &&\\
%$p=50$ & FP & 0.0123& 0.1138 &&\\
%$\sigma_y=0.5$ & PE & \textbf{1.0986}& 1.3907& 1.8965& 3.5062\\
% & SE PE & 0.0148& 0.0156& 0.0110& 0.0309\\
%\hline
%$n=100$ & FN & 0.0000 & 0.0000 &&\\
%$p=120$ & FP & 0.0035& 0.0697 &&\\
%$\sigma_y=0.5$ & Mean PE & \textbf{1.2250}& 1.5164 &&\\
% & SE PE & 0.0164& 0.0159 &&\\
%\hline\hline
%\end{tabular}}%
%\end{center}
%\caption{\label{tabnonlin}Comparison of FAR to three alternative
%methods in three nonlinear simulation
% settings.}
%\end{table}

%t2 #&#
\begin{table}[b]
\tabcolsep=0pt
\caption{Comparison of FAR to three alternative methods in three
nonlinear simulation
settings}
\label{tabnonlin}
\begin{tabular*}{\tablewidth}{@{\extracolsep{\fill}}@{}lccccc@{}}
\hline
& & \multicolumn{1}{c}{\textbf{FAR}} & \multicolumn{1}{c}{\textbf{FPCA}} & \multicolumn{1}{c}{\textbf{ADD}} & \multicolumn{1}{c@{}}{\textbf{SIR}}\\
\hline
$n=100$ & FN & 0.0000& 0.0000 &&\\
$p=5$ & FP & 0.1833& 0.1300 &&\\
$\sigma_y=0.5$ & Mean PE & {0.9792}& 1.3108& 1.7408& \textbf{0.8688}\\
& SE PE & 0.0132& 0.0174& 0.0074& 0.0049
\\[6pt]
$n=100$ & FN & 0.0000 & 0.0000 &&\\
$p=50$ & FP & 0.0171& 0.1138 &&\\
$\sigma_y=0.5$ & PE & \textbf{1.1068}& 1.3907& 1.8965& 3.5062\\
& SE PE & 0.0164& 0.0156& 0.0110& 0.0309
\\[6pt]
$n=100$ & FN & 0.0000 & 0.0000 &&\\
$p=120$ & FP & 0.0064& 0.0697 &&\\
$\sigma_y=0.5$ & Mean PE & \textbf{1.2108}& 1.5164 &&\\
& SE PE & 0.0157& 0.0159 &&\\
\hline
\end{tabular*}
\end{table}

The simulation results are summarized in Table~\ref{tabnonlin}, with
bold font indicating the statistically best prediction errors. As with
the linear setting it was not computationally feasible to implement ADD
or SIR for dimensionality $p$ larger than the sample size $n$. In the
low-dimensional setting of $p=5$, SIR produced the lowest mean
prediction error with FAR the second best. For the higher-dimensional
setting of $p=50$, the mean prediction error of SIR increased
dramatically and was the largest among all competitors. In the last two
settings, FAR was significantly superior to all three competing
methods, with FPCA generally providing the next best results. However,
we remark that the FPCA method is significantly slower than FAR in
these nonlinear settings due to the extra tuning parameter.

%f1 #&#
\begin{figure}[b]

\includegraphics{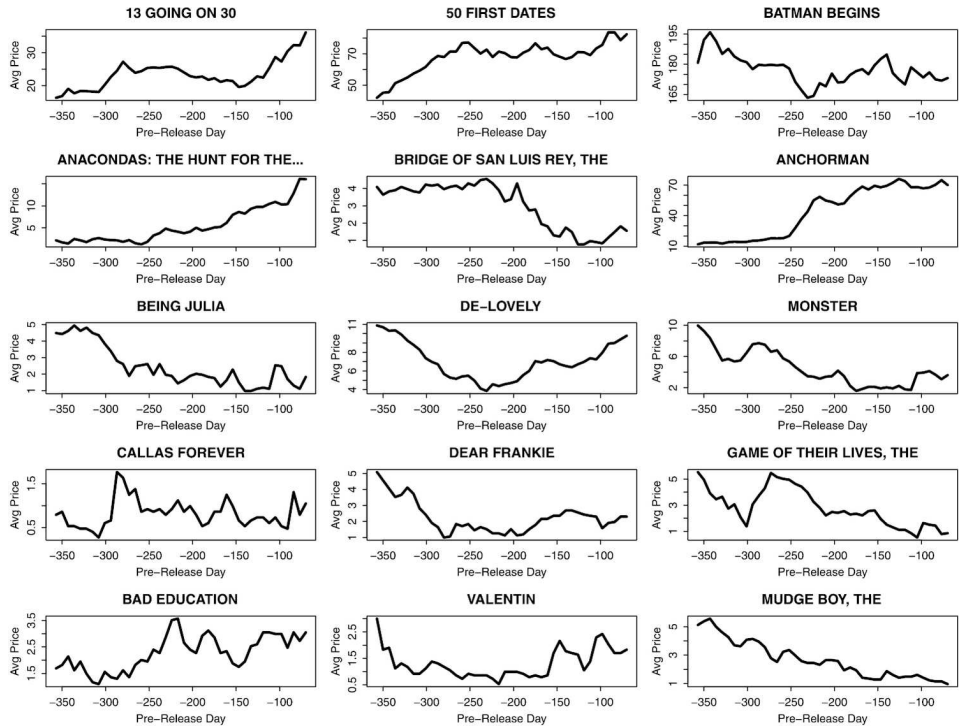}

\caption{Trading histories for a sample of movies from the HSX data
set.} \label{HSX}
\end{figure}

%In terms of prediction error, both FAR methods were superior to the
%ADD, MAVE and SIR approaches in all settings. In addition FAR was
%either superior or statistically tied with the FPCA approach in all
%simulations except the $\gamma=99\%$ situation. In this last setting
%FPCA provided some improvement in prediction accuracy, though no clear
%advantage in terms of model selection accuracy. By contrast FAR
%produced very large improvements in prediction accuracy in the $
%\gamma=69\%$ situation. The ADD, MAVE and SIR methods were competitive
%with FPCA in the $\gamma=69\%$ and $p=20$ simulation but were clearly
%inferior in the other settings. In the case of MAVE and SIR at least
%part of the explanation for their relatively worse performance is
%likely to be the fact that they do not perform any kind of variable
%selection.

%superior,
%though the improvement over FAR was relatively small. The linear FPCA
%methods were inferior to both FAR and FPCA.NL but outperformed the
%npenalized least squares approach. In comparing the shrunk and unshrunk
%versions of FAR, there was very little difference in the prediction
%errors. In general the unshrunk version resulted in lower false
%positive
%rates and, with the exception of the $n=100$ setting, similar false
%negative rates.

%s5 #&#
\section{Real data}
\label{realsec}

%In Sections~\ref{sechsx} and~\ref{meg} we respectively demonstrate
%FAR on medium and high dimensional data sets.

%s5.1 #&#
\subsection{Hollywood stock exchange data}
\label{sechsx}

%Our first data set comes from the {\it Hollywood Stock Exchange}
%(HSX), one of the best known online virtual stock markets. Online
%virtual stock markets operate in ways very similar to real life
%markets except that they are not necessarily based on real currency
%(i.e. participants often use
%virtual currency to make trades), and that each stock corresponds to
%an event or a parameter (rather than a company's shares). For
%instance, a value of 54 cents for the stock \emph{A democratic
%candidate will win
%the Presidential election} could be interpreted as the traders'
%collective belief that the democratic candidate has a 54\% chance of
%winning. If in fact the democratic candidate wins, then traders
%holding the
%democratic candidate's stock will liquidate (or ``cash-in") at \$1 per
%share; otherwise they receive \$0.

The goal for this analysis was to compare the accuracy of FAR and FPCA
in predicting the total box office revenue (over the first ten weeks
after release) for $262$ movies. We use pre-release trading histories
from the \textit{Hollywood Stock Exchange} (HSX), one of the best known
online virtual stock markets, as our functional predictors.
The Hollywood stock exchange has nearly 2 million active participants
worldwide. Each trader is initially endowed with \$2 million virtual
currency and can increase his or her net worth by strategically
selecting and trading movie stocks (i.e., buying low and selling high).
%Traders are further motivated by opportunities to exchange the accrued
%currency for merchandize and to appear on the daily {\it Leader Board}
%that features the most successful traders.
Figure~\ref{HSX} shows the HSX trading histories, between $52$ and
$10$ weeks prior to a movie's release, for a sample of $15$ out of the
$262$ movies in our data set. Each curve represents the traders'
collective daily average predictions of the box office revenue that the
movie will generate after it is released. In addition to the \textit{Daily
Average} curves, we also observed four additional predictors for each
movie: \textit{Accounts Trading}; \textit{Accounts Trading Short};
\textit{Shares Held
Short}; \textit{Shares Traded Sell}.

%f2 #&#
\begin{figure}[b]

\includegraphics{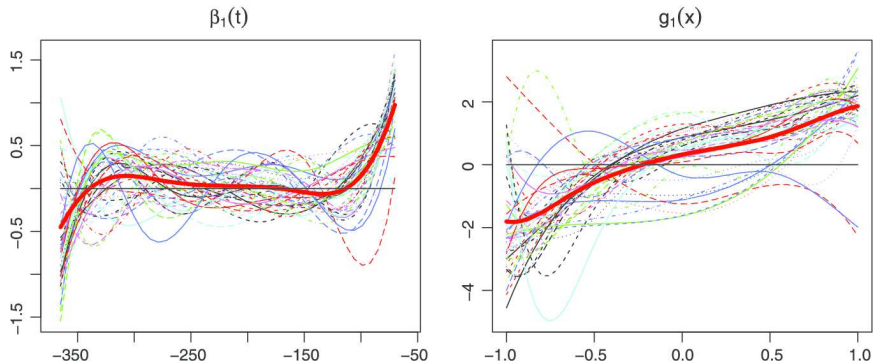}

\caption{The $\beta(t)$ and $g(x)$ curves corresponding to the Daily
Average variable in the HSX data.} \label{gbetaplot}
\end{figure}

%The goal for this analysis was to use FAR to predict the first ten
%weeks of cumulative box office revenue after release, based on the
%seven HSX trading history curves. Previous research has shown that
%online virtual stock markets such as HSX can provide impressive levels
%of predictive accuracy for not only movie returns but also in a
%variety of other settings \cite{foutz1}.
We only consider HSX curves from $10$ weeks prior to release date
because the goal is to form accurate revenue predictions early enough
to affect strategic decisions, such as, advertising budget, locations
of theater release, etc. We implemented the nonlinear versions of both
FAR and FPCA on the log revenues as this appeared to give superior
results for both methods. For FAR, we needed to select $3$ tuning
parameters, $\lambda$, $q$ and $d$, and for FPCA we also had $3$
tuning parameter, $\lambda$, $K$ and $d$. Hence, we randomly divided
the $262$ movies into three approximately equal partitions. The methods
were trained on the first group over grids of the tuning parameters,
the second group was used to select the final tuning parameters and out
of sample error rates were computed on the final group.

The mean hold out (log) prediction error, averaged over 50 random
partitions, was $2.45$ for FAR, while the FPCA error rate was higher at
$2.66$. The standard error in the difference between the FAR and FPCA
methods over the 50 random partitions was $0.10$. Both FAR and FPCA
chose Daily Average in all $50$ partitions, with the average model size
of FAR being $1.92$ and the average model size of FPCA being $1.86$.
The mean hold out (log) prediction error on the test movies using the
null model is $4.75$, indicating that using these functional predictors
from the trading histories indeed improves the prediction results.

%Table~\ref{tabhsx} displays the FAR and FPCA results, averaged over
%ten random partitions of the data. The first row corresponds to the
%average mean squared error on the test data. The standard shrunk
%version of FAR provided the best results, with the unshrunk FPCA
%method the closest competitor. The standard error in the difference
%between the FAR and FPCA* methods over the ten random partitions was
%approximately $0.11$. The remaining rows provide the proportion of
%times each variable was selected over the ten partitions. The Daily
%Average variable was selected in almost every model, while the other
%variables were selected sparingly. The standard shrunk FPCA method
%selected the largest models while FAR* produced the smallest models.

%\begin{table}[t]
% \centering
% \begin{tabular}{l\right\vert rrrr}
%& {\footnotesize FAR*} & {\footnotesize FAR} & {\footnotesize FPCA*} &
%{\footnotesize FPCA}\\\hline
%Mean PE & 2.771 & 2.437 & 2.602 & 2.667\\
%%SE PE & 0.287 & 0.171 & 0.169 & 0.164\\
%Daily Average & 1.0 & 1.0 & 0.9 & 1.0 \\
%Accounts Trading & 0.0 & 0.0 & 0.1 & 0.5\\
%Shares Traded & 0.1 & 0.2 & 0.0 & 0.8\\
%Shares Traded Cover & 0.1 & 0.2 & 0.0 & 0.6\\
%Dollar Volume Traded & 0.0 & 0.0 & 0.2 & 0.1\\
%Dollar Volume Traded Short & 0.0 & 0.2 & 0.0 & 0.1\\
%Dollar Volume Holding Short & 0.0 & 0.0 & 0.0 & 0.4
% \end{tabular}%
% \caption{Comparison of FAR to FPCA averaged over ten random
%partitions of the HSX movie data. We use * to denote unshrunk
%estimators.}
% \label{tabhsx}
%\end{table}%

Figure~\ref{gbetaplot} plots the 50 estimated $\beta(t)$ and $g(x)$
functions corresponding to the Daily Average variable with the solid
red lines representing the average effect. Most of the curves show
remarkably consistent patterns; $g(x)$ is estimated as a strictly
increasing, but nonlinear function, and $\beta(t)$ places
approximately zero weight on the earlier trading history and a larger
positive weight on roughly the final month under consideration. These
curves conform to our intuition that the trading history closest to
release date provides the strongest prediction accuracy and that there
is a positive correlation between HSX curves and movie revenues. The
nonlinear shape of $g(x)$ also suggests that a linear model would not
provide accurate results for this data.

%s5.2 #&#
\subsection{MEG data}
\label{meg}
Our second data set consisted of Magnetoencephalography (MEG)
recordings for $20$ subjects conducted at the Center for Clinical
Neurosciences, University of Texas Health
Science Center at Houston. The MEG readings for each subject were
recorded over $248$ ``channels'' at $356$ equally spaced time points.
Each channel measured the
intensity level of the magnetic field at a particular point on the
brain. Multiple trials, consisting of reading a patient a word and
measuring the MEG over time, were
recorded for each patient. We averaged the trials for each patient to
produce $248$ functional predictors, one for each channel. The response
of interest was whether the patient
was left ($14$ subjects) or right ($6$ subjects) brain dominated. We
coded $Y=1$ and $Y=-1$, respectively, for left- and right-brained
subjects. Some channels were missing
for some patients and were removed from the study, leaving a total of
$p=199$ predictors.

\begin{figure}

\includegraphics{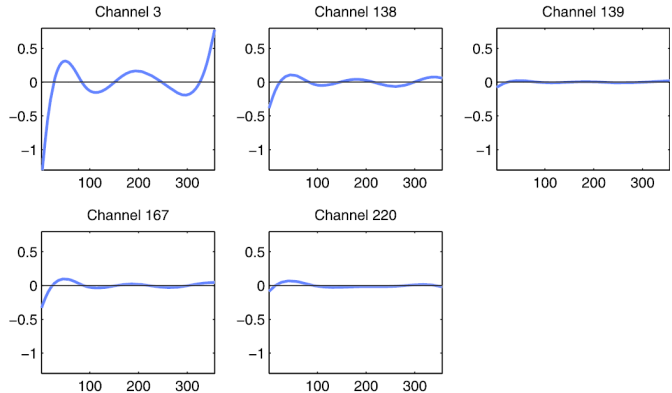}

\caption{Plots of $\beta(t)$ for linear FAR on the MEG data.}\label{figavglinear}
\end{figure}

This was a very challenging data set because the ratio of predictors to
observations was 10:1. We first fit the linear version of FAR to the
full data set using a five-dimensional basis for $\beta_j(t)$. The
tuning parameter, $\lambda$, was chosen as the point which minimized
the classification error using $20$-fold cross-validation. In this
setting, FAR selected
only a five variable model (Channels $3, 138, 139, 167$ and $220$),
which corresponded to a $20\%$ cross-validated error rate. Figure~\ref
{figavglinear} displays the
$\beta(t)$ curves for each selected channel. All five channels put the
bulk of their weight on the early time points. Channel $3$ appears to
provide the majority of the
predictive power with smaller contributions from Channels $138$ and
$167$. In particular $\beta_3$(t) represents a contrast between early
and late time points. Hence, people who start low in Channel 3 and end
high are predicted to be left-brained while the opposite is true for
right-brained patients.

%f3 #&#

We also fit the nonlinear version of FAR. Given the small number of
observations and the extra demands of fitting a nonlinear regression
method we felt it was prudent
to first perform a marginal pre-screening to select a smaller subset of
predictors for the final analysis. The marginal screening was performed
by running nonlinear FAR,
using a $7$-dimensional basis function, separately on each of the $194$
predictors that linear FAR did not choose and selecting the $45$ best
predictors in terms of
marginal prediction accuracy. Nonlinear FAR was then run on the $50$
predictors, including the $5$ selected by linear FAR. $20$-fold cross
validation was again used to
select the tuning parameter, resulting in five channels being selected.
The channels were not the same as those selected by linear FAR. The
cross-validated error rate was
$25\%$, suggesting that linear FAR may have a slight advantage on this data.

%s6 #&#
\section{Discussion}
\label{discsec}

FAR extends the recent linear penalized regression literature by
incorporating functional predictors and modeling general nonlinear
relationships. It has several
advantages over current functional regression methods. First, the
penalized approach automatically deals with high-dimensional data using
an efficient coordinate descent
algorithm. Second, the single index formulation provides a nonlinear
supervised method for projecting the predictors into a
lower-dimensional space, providing more
accurate results than the traditional linear unsupervised PCA approach.
Third, our theoretical results suggest that FAR should provide accurate
variable selection and prediction results and the
simulation results show that FAR outperforms traditional approaches.

There are three obvious possible extensions for FAR. The first is to
incorporate FAR into the generalized linear models setting.
Conceptually, such an extension could be
achieved by replacing the sum of squares term in (\ref{nonlinfar})
with the log likelihood and then using a modified version of the
coordinate descent algorithm to
maximize the criterion. The second possible extension would be to
replace the single index model with a multiple index model of the form,
$f_j(X_{ij}) =
\sum_{k=1}^Kg_{jk} (\int\beta_{jk}(t)X_{ij}(t)\,dt )$. This
would increase the flexibility of FAR to model more general nonlinear
relationships. Finally, FAR could be extended to model functional
responses in addition to functional predictors.

%\section*{Acknowledgements}
%We would like to thank the Associate Editor and two referees for many
%helpful suggestions that
%improved the paper.
%This work was partially supported by NSF Grant DMS-0906784.

%\section{Initial Parameters for the FAR Algorithm}
%\label{initialvalues}
%
%%We use the following coordinate descent algorithm to minimize
%%(\ref{etaopt}).\\
%
%To get initial values for the $\bfeta_j$ parameters, we first pretend
%that $g_j(t)$'s were linear functions and use the linear FAR algorithm
%to get the estimates $\hbfeta_j$'s. Then, for each~$j\in\{1,...,p\}$,
%if $\left\Vert\hbfeta_j\right\Vert_2 = 0$, we update $\hbfeta_j$ to
%be the loading
%vector of the first principle component decomposition of $\Theta_j$.
%%\begin{enumerate}
%%\item Fix all $\hf_k$ for $k\ne j$ and compute
%% the residual vector $\bR_j = \bY- \sum_{k\ne j} \hf_k(X_k)$.
%%\item Use Sliced Inverse Regression \cite{li1}, with $\btheta_{ij}$
%as the
%% predictor, to estimate $\bfeta_j$.
%%\item Given $\hat{\bfeta}_j$, estimate $\bxi_j$ using the usual
%least
%% squares solution to (\ref{partialetaopt}).
%%\end{enumerate}
%%After Step 3. the FAR algorithm operates in an identical fashion for
%all
%%values of $\lambda$.

\begin{appendix}
\section{Details of the nonlinear FAR algorithm}\label{stepA}
In the initialization step (step~0) of this algorithm, some of the
$\bfeta_j$'s will likely be set to zero. This suggests that the
corresponding predictors do not appear related to the response.
However, the initialization assumes a linear model. It is conceivable
that a response that appears unimportant using a linear model will
become statistically significant using a nonlinear model. Hence, if
$\bfeta_j$ is estimated to be zero in step~0 we instead set $\bfeta
_j$ equal to the loading vector of the first principal component of
$\Theta_j$. This estimate is the direction that explains the most
variability in $X_{ij}(t)$ so is the most natural unsupervised
projection and allows for potential nonlinear relationships to be
detected in step~2.

To implement step~3 of the FAR algorithm, we minimize (\ref{etaopt})
with respect to the $\bfeta_j$'s. Directly minimizing (\ref{etaopt})
is difficult due to the nonlinearity of the functions $g_j(t)\approx
\bh(t)^T\bxi_j$. To overcome this difficulty, we observe that, with
the estimate $\hbxi_j$ from step~2 and the current value $\bfeta
_{j,\mathrm{old}}$ of $\bfeta_j$, the first-order approximation of $g(\btheta
_{ij}^T\bfeta_j)\approx\bh(\btheta_{ij}^T\bfeta_j)^T\hbxi_j$ is
%
%e16 #&#
\begin{equation}
\label{eqhapprox} \qquad\bh\bigl(\btheta_{ij}^T\bfeta_j
\bigr)^T\hbxi_j\approx\bh\bigl(\btheta
_{ij}^T\bfeta_{j,\mathrm{old}}\bigr)^T
\hbxi_j + \bh'\bigl(\btheta_{ij}^T
\bfeta_{j,\mathrm{old}}\bigr)^T\hbxi_j \cdot
\btheta_{ij}^T(\bfeta_j-\bfeta_{j,\mathrm{old}}).
\end{equation}
Thus, we can approximate (\ref{etaopt}) as
%
%e17 #&#
\begin{equation}
\label{etaoptapprox} \sum_{i=1}^n
\Biggl(R_i - \sum_{j=1}^p
\bh'\bigl(\btheta_{ij}^T\bfeta_{j,\mathrm{old}}
\bigr)^T\hbxi_j \cdot\btheta_{ij}^T(
\bfeta_j-\bfeta_{j,\mathrm{old}}) \Biggr)^2,
\end{equation}
where\vspace*{1pt} $R_i = Y_i - \sum_{j=1}^p \bh(\btheta_{ij}^T\bfeta
_{j,\mathrm{old}})^T\hbxi_j $, that is, the residual for the $i$th observation
from step~2 of the algorithm in the current iteration. The above
approximation (\ref{etaoptapprox}) is a quadratic function of $\bfeta
_j$ and can be minimized easily. Hence, the new value of $\bfeta_j$ is
updated as the minimizer of (\ref{etaoptapprox}). We also note that if
the estimate $\hbxi_j$ from step~2 is $\bzero$, then the
corresponding value of $\bfeta_j$ will not be updated.

%s8 #&#
\section{Technical conditions of Theorems \texorpdfstring{\protect\ref{T2}--\protect
\ref{T3}}{1--2}}\label{linearCon}
\label{appb}
%\subsection{Technical Conditions}

We make the following assumption on the functional predictors
$X_{ij}(t)$ and the corresponding regression coefficients $\beta_j(t)$.

%co1 #&#
\begin{assumption}\label{as3}
\textup{(A)} Functional predictors, $\{X_{ij}:[0,1]\rightarrow\RR,
i=1,\ldots,n, \break j=1,\ldots,{p_n}\}$, belong to a Sobolev ellipsoid of
order two: there exists a universal constant $C$, such that $\sum
_{k=1}^{\infty} \theta^2_{ijk}k^4\le C^2$ for all $i=1,\ldots,n,
j=1,\ldots,{p_n}$.

(B) The true coefficient functions satisfy $\max_{ j\in\Mfrak
_0}\int_0^1\beta^2_j(t)\,dt \leq\widetilde C$ with $\widetilde C$ some
positive constant.
\end{assumption}

%Thus, it follows from Condition~\ref{as3}(B) that uniformly for all
%true predictor $j\in\Mfrak_0$,
%\begin{eqnarray*}
%\int_0^1\beta_j(t)^2\,dt = \sum_{l=1}^\infty\eta_{0,jl}^2 \leq\widetilde C,
%\end{eqnarray*}
%while for each noise predictor $j\notin\Mfrak_0$, $\eta_{0,jl}=0$ for
%all $l=1,2,\cdots$.
%Using \eqref{eqbasis-rep}, the $j$th additive component has the
%following representation
%\begin{eqnarray}\label{fexpan}
%f_j(X_{ij}) = \int_0^1 X_{ij}(t)\beta_j(t)\,dt = \sum_{l=1}^\infty
%\theta_{ijl}\eta_{0,jl}.
%\end{eqnarray}
%For a given sequence of integers $q_n$ depending only the sample size
%$n$, write $\bfeta_{0j} = (\eta_{0,j1},\cdots,\eta_{0,jq_n})\t$.
%Correspondingly, write $\btheta_{ij} = (\theta_{ij1},\cdots,
%\theta_{ijq_n})\t$. Thus, the $j$th additive component $f_j(X_{ij})$
%can be approximately as $\btheta_{ij}\t\bfeta_{0j}$. Denote by
%$e_{ij}$ the approximation error, that is,
%\begin{eqnarray*}
%e_{ij} = f_j(X_{ij})-\btheta_{ij}^T\bfeta_j.
%\end{eqnarray*}
%Then uniformly across all $i$ and $j\in\Mfrak_0$,
%\begin{eqnarray}
%\nonumber\right\vert e_{ij}\right\vert^2& = \Big\right\vert
%\sum_{l=q_n+1}^\infty\theta_{ijl}l^{-2}l^2
%\eta_{0,jl}\Big\right\vert^2\leq\sum_{l=q_n+1}^\infty
%\eta_{0,jl}^2l^{-4}
%\sum_{l=q_n+1}^\infty\theta_{ijl}^2l^{4} \\
%& \leq& C^2 q_n^{-4}\sum_{l=q_n+1}^\infty\eta_{0,jl}^2 \leq\widetilde C
%C^2q_n^{-4}, \label{eqapp-err}
%\end{eqnarray}
%where the last two inequalities are because of Condition~\ref{as3}(A)
%and (B), respectively.
%

Note that the linear FAR model can be written as
%
%e18 #&#
\begin{equation}
\label{linearapprox} Y_i = \sum_{j=1}^p
\Theta_{j}\bfeta_j + \veps_i^*,
\end{equation}
where $\veps^*_i = \veps_i + \sum_{j=1}^pe_{ij}$ with $e_{ij}$
defined in (\ref{deferr}).
When $j\in\Mfrak_0^c$, $\beta_j(t)=0$ and thus the approximation
error $e_{ij}$ in (\ref{fexpan}) disappears. Thus, in view of (\ref
{eqapp-err}), the approximation error satisfies that
\begin{eqnarray*}
&&\Biggl\llvert\sum_{j=1}^pe_{ij}
\Biggr\rrvert\leq\sum_{j \in\Mfrak_0}\llvert e_{ij}
\rrvert\leq Cs_nq_n^{-2},
\end{eqnarray*}
uniformly over all $i = 1, \ldots, n$.

%Condition~\ref{as3} puts restrictions on the number of nonzero
%parameters $s_n\times q$, that is, $s_n \times q = o(n)$ and
%$s_n=o(n^{2\delta-\frac{1}{2}})$ for some
%$\delta> \frac{1}{4}$. This condition is imposed because, in order to
%ensure that the approximation holds uniformly, the number of true
%predictors, $s_n$, cannot grow
%too fast, and to avoid over-fitting, the number of basis functions,
%$q$, must also be constrained. It is easy to derive from Condition~
%\ref{as3} that $e_{ij}$, defined in
%(\ref{fexpan}), satisfies $\max_{i,j\in\Mfrak_0}\left\vert e_{ij}\right
%\vert= o(n^{-2
%\delta})$ and $e_{ij}=0$ for $j\in\Mfrak_0^c$. Hence, since it is
%assumed that
%$s_n=o(n^{2\delta-\frac{1}{2}})$, then the mean of the error term $
%\veps_i^*$ in (\ref{linearapprox}) has order $o(n^{-1/2})$ uniformly
%across $i$. %If $\bb(t)$ is the true basis from which $X_{ij}(t)$ and
%$\beta_j(t)$ are generated then the approximation errors are $0$ and
%Condition~\ref{as3} is not needed. However, in general assuming a
%perfect representation seems to be a strong assumption.

Our second set of conditions concern the shape of the penalty function,
the strength of the signal and the correlation structure of the predictors.
%We next show that with probability tending to 1, there exists a local
%minimizer of %(\ref{e005}) whose support is $\Mfrak_0$. We need the
%following assumption:

%co2 #&#
\begin{assumption}\label{as2}
(A) For any fixed $\lambda>0$, $\rho_{\lambda}(t)$ is
concave and nondecreasing in $[0, \infty)$, and has nonincreasing
first derivative
$\rho_{\lambda}'(t)$. Further, $\rho'_\lambda(0+)>0$.

(B) Let $a_n=\min_{j\in\Mfrak_0}\llVert \Theta_j\bfeta
_{0,j}\rrVert
/\sqrt{n}$. It holds that $n^{\alpha} a_n \rightarrow\infty$ with
$\alpha\in(0,\frac{1}{2})$.

(C) It holds that $\rho'_{\lambda_n}(a_n/2)=o (n^{-\alpha
}q_n^{-1}s_n^{-1/2} )$ and $\sup_{t\geq\sfrac{a_n}{2}}\rho
''_{\lambda_n}(t)=o(1)$.

(D) There exists a positive constant $c_0$ such that
%
%e19 #&#
\begin{equation}
\label{eqeigen} c_0\leq\min_{j\in\Mfrak_0}
\Lambda_{\min}\biggl(\frac{1}{n} \Theta_{j}^T
\Theta_{j}\biggr)<\Lambda_{\max}\biggl(\frac{1}{n}
\Theta_{\Mfrak_0}^T\Theta_{\Mfrak_0}\biggr)\leq
c_0^{-1},
\end{equation}
where $\Lambda_{\min}$ and $\Lambda_{\max}$
are the smallest and largest eigenvalues of a matrix, respectively.
Further, we have
%
%e20 #&#
\begin{eqnarray}
\label{e015}
&& \max_{j\in\Mfrak_0^c} \bigl\llVert\Theta_j
\bigl(\Theta_j^T\Theta_j
\bigr)^{-1}\Theta_j^T \Theta_{\Mfrak_0}
\bigl(\Theta_{\Mfrak_0}^T\Theta_{\Mfrak_0}
\bigr)^{-1} \bigr\rrVert_{\infty,
2}<\frac{\sqrt{c_0}}{2\sqrt{n}}
\frac{\rho'_{\lambda_n}(0+)}{\rho
'_{\lambda_n}(a_n/2)},%\\
%\label{e016}& \max_{j\in\Mfrak_0^c}\left\Vert\Theta_j(\Theta_j^T
%\Theta_j)^{-1}
%\right\Vert_{\infty, 2}< \frac{\rho'_{\lambda_n}(0+)}{2\sqrt{n}
%\lambda_n},
\end{eqnarray}
where for a matrix $B$, $\llVert B\rrVert _{\infty,2} = \sup
_{\llVert {\mathbf x}\rrVert
_\infty=1}\llVert B{\mathbf x}\rrVert _2$ with ${\mathbf x}$ a vector.

(E) The model errors $\veps_i$, $i=1,\ldots, n$ are
independent and identically distributed as $N(0,\sigma^2)$.
\end{assumption}

Condition~\ref{as2}(A) requires that the penalty functional, $\rho
_{\lambda}(t)$, is concave and singular at 0. Many penalty functions
proposed in the literature such as the hard thresholding penalty, SCAD
\cite{fan2} and SICA \cite{lv1} all satisfy this condition. From
(\ref{fexpan}), we see that Condition~\ref{as2}(B)
places a lower bound on the signal strength of the true predictors
$j\in\Mfrak_0$. In particular, it assumes that the weakest signal,
$a_n$, can decay with sample
size but the decay rate cannot be faster than $n^{-\alpha}$. Condition
\ref{as2}(C) is a mild condition which can be easily satisfied by
penalty functions with flat tails. For instance, if $\lambda_n =
o(a_n/2)$, then for SCAD penalty, it can be verified from the
definition that $\rho'_{\lambda_n}(a_n/2) = 0$ and $\rho''_{\lambda
_n}(t) = 0$ for all $t \geq a_n/2$, and thus\vspace*{1pt} Condition~\ref{as2}(C) is
satisfied. Although Condition~\ref{as2}(C) assumes the existence of
the second-order derivative for
$\rho_{\lambda_n}(t)$, it can be relaxed to the existence of the
first-order derivative by using the local concavity definition in \cite
{lv1}. Condition~\ref{as2}(D)
relates to the design matrix for the signal predictors, $\Theta
_{\Mfrak_0}$. We assume that the eigenvalues for the design matrix
corresponding to true predictors
are bounded from below and above. If $\Theta_{\Mfrak_0}$ is
orthogonal, then (\ref{eqeigen}) is satisfied with $c_0 = 1$. The
upper bound in condition (\ref{e015}) depends on the penalty function
through the ratio
$\rho'_{\lambda_n}(0+)/\rho_{\lambda_n}'(a_n/2)$, which is larger
than 1 for concave penalties and equal to 1 for the group Lasso penalty,
$\rho_{\lambda_n}(t)=\lambda_n t$. For instance, if $\lambda_n =
o(a_n)$, then $\rho'_{\lambda_n}(0+)/\rho_{\lambda_n}'(a_n/2) =
\infty$ for SCAD penalty and thus (\ref{e015}) is satisfied
automatically. The detailed proofs of Theorems~\ref{T2} and~\ref{T3}
are in the supplementary materials \cite{FAR}.

%s9 #&#
\section{Technical conditions and proof of Theorems~\texorpdfstring{\lowercase{\protect\ref{nonlinthm}}--\lowercase{\protect\ref{threshthm}}}{3--4}}
\label{appnonlimthm}

%s9.1 #&#
\subsection{Conditions}

Given an orthonormal basis expansion for $\beta_j(t)$, that is, $\beta
_j(t)=\sum_{l=1}^{\infty}\eta^*_{jl}b_l(t)$, we will define $\bfeta
_j^*=(\eta^*_{j1},\ldots,\eta^*_{jq_n})^T$.
%Recall that we write $\btheta_{ij}$ for $(\theta_{ij1},...
%\theta_{ijd_n})^T$, using a similar expansion for~$X_{ij}(t)$.
We will also define $f_j^*(\btheta_{ij})=\bh_{\bfeta^*_j,j}(\btheta
_{ij}^T\bfeta_j^*)\bxi^*$, where $\bxi_j^*$ is chosen to minimize
$\sum_{i=1}^n [\bh_{\bfeta^*_j,j}\times\break (\btheta_{ij}^T\bfeta_j^*)\bxi
-g_{0j}(\btheta_{ij}^T\bfeta_j^*)]^2$ over~$\bxi\in\RR^{d_n}$ with
the constraint $\sum_i f_j^*(\btheta_{ij})=0$. Note that~$f_j^*$,
$\bxi_j^*$, $\bfeta_j^*$ and $\btheta_{ij}$ depend\vspace*{1pt} on~$n$, but we
omit the corresponding subscripts for the simplicity of the notation.
The following are the technical conditions for the theory in
Section~\ref{nonlintheorysec}. A discussion of the conditions is
given below.

%co3 #&#
\begin{assumption}\label{cond3}
(A) Functional predictors, $\{X_{ij}:[0,1]\rightarrow\RR,
i=1,\ldots,n, j=1,\ldots,{p_n}\}$, belong\vspace*{1pt} to a Sobolev ellipsoid of
order two: there exists a universal constant $C$, such that $\sum
_{k=1}^{\infty} \theta^2_{ijk}k^4\le C^2$ for all $i=1,\ldots,n,
j=1,\ldots,{p_n}$.

(B) The true index functions, $\{\beta_j(t), j\in\Mfrak_0\}
$, satisfy $\int_0^1\beta^2_j(t)\,dt=1$.

(C) Errors~$\varepsilon_i$ are independent and uniformly
sub-Gaussian.

(D) The true link functions, $g_{0j}$, are twice continuously
differentiable and are bounded, together with their first and second
derivatives, uniformly over~$j\in\Mfrak_0$ and~$n$.

(E) For\vspace*{1pt} each $\bfeta$ with $\llVert \bfeta\rrVert =1$ and each
$j\le
{p_n}$ let~$Q_{\bfeta,j,n}$ denote the empirical distribution
associated with the index values $\bfeta^T\btheta_{1j},\ldots,
\bfeta^T\btheta_{nj}$. Assume that there exist corresponding
probability distributions~$P_{\bfeta,j,n}$, each with bounded support
and a positive continuous density, such that the densities are bounded
both above and away from zero uniformly over $j$ and $n$, and
%
%e21 #&#
\begin{equation}
\sup_{u\in\RR, \llVert \bfeta\rrVert =1, 1\le j\le{p_n}}\bigl\llvert
Q_{\bfeta
,j,n}(-
\infty,u]-P_{\bfeta,j,n}(-\infty,u]\bigr\rrvert=o \bigl(d_n^{-1}
\bigr).
\end{equation}
\end{assumption}

Condition~\ref{cond3}(A) is identical to Condition~\ref{as3}(A),
imposed for the linear FAR theory. It is a common smoothness
requirement in nonparametric regression, when the orthogonal basis
approach is used, as discussed, for example, in Chapter~8 in \cite
{Bwasserman1}. Condition~\ref{cond3}(B) is imposed for
identifiability. Conditions~\ref{cond3}(C) and~(D) are typical in
high-dimensional regression and nonparametric regression problems,
respectively. The reason we require uniformity is to handle the
situation where the number of signal predictors grows with~$n$. Again,
uniformity is needed to handle the growing number of signal predictors.\vspace*{1pt}
Condition~\ref{cond3}(E) ensures that the candidate index values,
$\bfeta^T\btheta_{1j},\ldots, \bfeta^T\btheta_{nj}$, have
sufficiently regular distributions. Assumptions of this form are
typical in spline estimation \cite{zhou2}, for example.

We impose two more assumptions below. Condition~\ref{as4}(A) is a
natural generalization of the \emph{compatibility condition} used in
high-dimensional additive models, for example, in \cite{meier1} and
Section~8.4. in \cite{Bbuhlmann1}. Note that because we do not use a
smoothness penalty in our estimation approach, the smoothness penalty
does not appear in the compatibility condition. Condition~\ref{as4}(B)
is a version of the standard regularity condition on the behavior of
the sum of squares function near its minimum. Assumptions of this form
have been imposed in the single index model literature, for example,
\cite{yu1}. We again require uniformity over $j\in\Mfrak_0$ to
handle the growing number of signal predictors.

%co4 #&#
\begin{assumption}\label{as4}
(A) There exists a positive universal constant $\phi^2$ for
which the following holds. If functions $\{h_j, j=1,\ldots,{p_n}\}$
are such that each $h_j$ is a difference of two functions in $\F_j^0$,
and inequality $\sum_{j\in\Mfrak_0^c}\llVert h_{j}\rrVert _n\le3\sum
_{j\in
\Mfrak_0}\llVert h_{j}\rrVert _n$ is satisfied, then the following inequality
holds: $\sum_{j\in\Mfrak_0}\llVert h_{j}\rrVert _n^2\le\llVert \sum
_{j=1}^{{p_n}}h_j\rrVert ^2_n/\phi^2$.

(B) There exist positive universal constants $\tau$, $c_1$ and
$c_2$, such that for all sufficiently large~$n$ and each $f_{\bfeta
_j,j}(\cdot)=\bh_{\bfeta_j,j}(\bfeta_j^T\cdot)\bxi$ with~$\llVert
\bfeta_j\rrVert =\llVert \bfeta_j^*\rrVert $ and $j\in\Mfrak_0$,
inequalities $\llVert
f_{\bfeta_j,j} - f^*_j\rrVert _n\le\tau$ and $\llVert \bfeta_j-\bfeta
_j^*\rrVert <\llVert
\bfeta_j+\bfeta_j^*\rrVert $ imply $\llVert \bfeta_j-\bfeta
_j^*\rrVert \le c_1\llVert
f_{\bfeta_j,j} - f^*_j\rrVert _n$ and $\llVert f_{\bfeta_j^*,j} -
f^*_j\rrVert _n\le
c_2\llVert f_{\bfeta_j,j} - f^*_j\rrVert _n$.
%\item[{\mathbf A7}.] If we let $({\mathbf a},\D)=(d_n^{-1/2}(\bxi-
%\bxi^*)^T,(
%\bfeta-\bfeta^*)^T)$ and define function $\kappa_{j}({\mathbf a},\D)$
%as
%$(1/n)\sum_{i=1}^n \left[\bh_{\bfeta_j,j}(\bfeta_j^T\btheta_{ij})\bxi
%- f^*_j(\btheta_{ij})\right]^2$, then the eigenvalues of $
%\kappa_{j}''(0)$ are bounded away from zero uniformly over all
%sufficiently large~$n$ and $j\in\Mfrak_0$.
\end{assumption}

%s9.2 #&#
\subsection{Preliminaries}

% Recall that we treat predictors as deterministic.

We start by deriving a bound on the error due to our approximation of
index functions~$\beta_j$ and link functions~$g_{0j}$. Observe that
\begin{eqnarray*}
\biggl\llvert\int_0^1\beta_j(t)X_{ij}(t)\,dt
- \btheta_{ij}^T\bfeta_j^*\biggr\rrvert
^2&=&\Biggl\llvert\sum_{k=q_n+1}^{\infty}
\eta^*_{jk}k^{-2}k^2\theta_{ijk}\Biggr
\rrvert^2
\\
&\le&\sum_{k=q_n+1}^{\infty}\bigl(
\eta^*_{jk}\bigr)^2k^{-4}\sum
_{k=q_n+1}^{\infty}\theta^2_{ijk}k^{-4}.
\end{eqnarray*}
Condition~\ref{cond3}(A) implies that the right-most sum is bounded by
a universal constant~$C^2$. Also\vspace*{1pt} note that $\sum_{l=q_n+1}^{\infty
}(\eta_{jl}^*)^2k^{-4}\le q_n^{-4}\sum_{k=q_n+1}^{\infty}(\eta
_{jk}^*)^2\le q_n^{-4}$, by Condition~\ref{cond3}(B). Thus, if we set
$I_{ij}=\int_0^1\beta_j(t)X_{ij}(t)\,dt$, then the bound $\llvert I_{ij} -
\btheta_{ij}^T\bfeta_j^*\rrvert \le Cq_n^{-2}$ holds for all $n$, $i$, and
$j\in\Mfrak_0$.
Hence, if we let $\widetilde C$ be the uniform bound over the first
derivatives in Condition~\ref{cond3}(D), then
%
%e22 #&#
\begin{equation}
\label{prelimappr1} \bigl\llvert g_{0j}(I_{ij}) -
g_{0j} \bigl(\btheta_{ij}^T\bfeta_j^*
\bigr)\bigr\rrvert\le\widetilde C\bigl\llvert I_{ij} -
\btheta_{ij}^T\bfeta_j^*\bigr\rrvert= O
\bigl(q_n^{-2}\bigr),
\end{equation}
uniformly over $i$ and $j\in\Mfrak_0$. Set $M_j=\sup_t\llvert g^{\prime
\prime}_{0j}(t)\rrvert $ for $j\in\Mfrak_0$, and note that constants $M_j$
are uniformly bounded by Condition~\ref{cond3}(D). Taking advantage of
the approximation bounds for the cubic B-splines (e.g., Corollary~6.21
in \cite{BSchumaker1}), we then have
%
%e23 #&#
\begin{equation}
\label{prelimappr2} n^{-1}\sum_{i=1}^n
\bigl(g_{0j} \bigl(\btheta_{ij}^T
\bfeta_j^* \bigr)-f^*_j(\btheta_{ij})
\bigr)^2=O\bigl(d_n^{-4}M_j^2
\bigr)=O\bigl(d_n^{-4}\bigr),
\end{equation}
uniformly over $j\in\Mfrak_0$. Combining inequalities~(\ref
{prelimappr1}) and~(\ref{prelimappr2}), we deduce
$\llVert f_{0j}-f_{j}^*\rrVert _n=O(q_n^{-2}+d_n^{-2})=O(d_n^{-2})$,
uniformly over $j$.
Note that for $j\in\Mfrak_0^c$, both $f_{0j}$ and $f_{j}^*$ are zero.
Consequently,
%
%e24 #&#
\begin{equation}
\label{approxbnd} \sum_{j=1}^{{p_n}}\bigl\llVert
f_j^*-f_{0j}\bigr\rrVert_n=O
\bigl(s_nd_n^{-2}\bigr).
\end{equation}
This gives us a useful bound on the approximation error.

We will write~$f_0(\bX_i)$ for $\sum_{j=1}^{{p_n}} f_{0j}(X_{ij})$;
we also write~$\hat f(\btheta_i)$ for $\sum_{j=1}^{{p_n}}\hat f_j(\btheta_{ij})$ and define $f^*$ by analogy. To be consistent with
the standard least-squares estimation notation, we will
write~$(\varepsilon,f)_n$ for $n^{-1}\sum_{i=1}^n\varepsilon_i
f(\btheta_{i})$. We will need the following result, which is proved in
the supplementary material \cite{FAR}.

%
%le1 #&#
\begin{lemma}
\label{lementrcons}
Define $r_n=n^{-1/2}q_n^{1/2}+n^{-1/2}\sqrt{\log{p_n}}$.
There exists a positive universal constant $C_1$, such that
%
%e25 #&#
\begin{equation}
\bigl(\varepsilon,\hat f-f^*\bigr)_n\le C_1s_nr_n^2+C_1r_n
\sum_{j=1}^{p_n}\bigl\llVert\hat f_j - f_j^*\bigr\rrVert_n,
\end{equation}
with probability tending to one.
\end{lemma}

%s9.3 #&#
\subsection{Main body of the proof}

Let $\llVert y-f\rrVert ^2_n$ denote $n^{-1}\sum_{i=1}^n(Y_i-f(\btheta
_{i}))^2$
and let $\llVert f\rrVert ^2_n$ denote $n^{-1}\sum_{i=1}^n f(\btheta_{i})^2$.
Consider the following simple identity:
%
%e26 #&#
\begin{equation}
\label{basicineq1} \llVert y-\hat f\rrVert_n^2-\bigl
\llVert y-f^*\bigr\rrVert_n^2=\llVert\hat f-f_0\rrVert_n^2-\bigl\llVert
f^*-f_0\bigr\rrVert_n^2-2\bigl(\varepsilon,
\hat f - f^*\bigr)_n.
\end{equation}
Note that $\llVert y-\hat f\rrVert _n^2+\lambda_n\sum
_{j=1}^{p_n}\llVert \hat f_j\rrVert _n-\llVert y-f^*\rrVert _n^2-\lambda_n\sum
_{j=1}^{p_n}\llVert f_j^*\rrVert _n\le0$ by
the definition of~$\hat f$. Let~$e_n$ denote the approximation
error, $\llVert f^* - f_0\rrVert _n$. Inequality~(\ref{basicineq1})
then implies
\[
\llVert\hat f-f_0\rrVert_n^2+
\lambda_n\sum_{j=1}^{p_n}\llVert
\hat f_j\rrVert_n\le e_n^2+2
\bigl(\varepsilon,\hat f - f^*\bigr)_n+\lambda_n\sum
_{j=1}^{p_n}\bigl\llVert f_j^*
\bigr\rrVert_n.
\]
By Lemma~\ref{lementrcons}, the above inequality yields
%
%e27 #&#
\begin{eqnarray}\label{basicineq1b}
&& \llVert\hat f-f_0\rrVert_n^2+
\lambda_n\sum_{j=1}^{p_n}\llVert
\hat f_j\rrVert_n
\nonumber\\[-10pt]\\[-10pt]\nonumber
&&\qquad \le e_n^2+2C_1s_nr_n^2+2C_1r_n
\sum_{j=1}^{p_n}\bigl\llVert\hat f_j - f_j^*\bigr\rrVert_n+
\lambda_n\sum_{j=1}^{p_n}\bigl
\llVert f_j^*\bigr\rrVert_n,
\end{eqnarray}
with probability tending to one.

\begin{longlist}
\item[\textit{Case} (i).]
Consider the event $e_n^2+C_1s_nr_n^2\ge r_n\sum_{j=1}^{p_n} \llVert
\hat f_j-f_j^*\rrVert _n$.

Note that $e_n^2=O(s_n^2d_n^{-4})$ by~(\ref{approxbnd}).
Thus, $\sum_{j=1}^{p_n} \llVert \hat f_j-f_j^*\rrVert _n =
O(s_n^2n^{1/2}d_n^{-4}\*q_n^{-1/2}+s_nr_n)$. Consequently,
$\sum_{j=1}^{p_n} \llVert \hat f_j-f_{0j}\rrVert
_n=O(d_n^{-4}q_n^{-1/2}+s_nr_n+s_nd_n^{-2})$, which implies the
stochastic bound in display~(\ref{nonlinthmbnd}).\vspace*{1pt}

\item[\textit{Case} (ii).]
Consider the event $e_n^2+C_1s_nr_n^2< r_n\sum_{j=1}^{p_n} \llVert
\hat f_j-f_j^*\rrVert _n$.\vadjust{\goodbreak}

Using inequality $\llVert \hat f-f^*\rrVert ^2_n\le
2\llVert \hat f-f_0\rrVert _n^2 +2e_n^2$ together with~(\ref{basicineq1b}), we get
\begin{eqnarray*}
&& \bigl\llVert\hat f-f^*\bigr\rrVert_n^2+2
\lambda_n\sum_{j=1}^{p_n} \llVert
\hat f_j\rrVert_n
\\
&&\qquad \le4e_n^2+4C_1s_nr_n^2+4C_1r_n
\sum_{j=1}^{p_n} \bigl\llVert\hat f_j - f_j^*\bigr\rrVert_n+2
\lambda_n\sum_{j=1}^{p_n} \bigl
\llVert f_j^*\bigr\rrVert_n.
\end{eqnarray*}
On the event $e_n^2+C_1s_nr_n^2< r_n\sum_{j=1}^{p_n} \llVert \hat f_j-f_j^*\rrVert _n$ the above inequality simplifies to
%
%e28 #&#
\begin{eqnarray}\label{basicineq3}
&& \bigl\llVert\hat f-f^*\bigr\rrVert_n^2+2
\lambda_n\sum_{j=1}^{p_n} \llVert
\hat f_j\rrVert_n
\nonumber\\[-8pt]\\[-8pt]\nonumber
&&\qquad \le4(C_1+1)r_n
\sum_{j=1}^{p_n} \bigl\llVert\hat f_j - f_j^*\bigr\rrVert_n+2\lambda
_n\sum_{j=1}^{p_n} \bigl\llVert
f_j^*\bigr\rrVert_n.
\end{eqnarray}
Because we assume $r_n=O(\lambda_n)$, we can rewrite inequality~(\ref
{basicineq3}) as
%
%e29 #&#
\begin{equation}
\label{stochbndlam} \bigl\llVert\hat f-f^*\bigr\rrVert_n^2
= \sum_{j=1}^{p_n} \bigl\llVert\hat f_j-f^*_j\bigr\rrVert_n O(
\lambda_n).
\end{equation}
Inequality~(\ref{basicineq3}) also gives
\[
\sum_{j\in\Mfrak_0^c} \llVert\hat f_j\rrVert
_n\le2\lambda_n^{-1}(C_1+1)r_n
\sum_{j=1}^{p_n} \bigl\llVert\hat f_j - f_j^*\bigr\rrVert_n+\sum
_{j\in\Mfrak_0} \bigl\llVert\hat f_j-f_j^*
\bigr\rrVert_n.
\]
Consequently,
\begin{eqnarray*}
\sum_{j\in\Mfrak_0^c} \bigl\llVert\hat f_j-f_j^*\bigr\rrVert_n
&\le& 2\lambda
_n^{-1}(C_1+1)r_n\sum
_{j\in\Mfrak_0^c} \bigl\llVert\hat f_j-f_j^*
\bigr\rrVert_n
\\
&&{} +\bigl[2\lambda_n^{-1}(C_1+1)r_n+1
\bigr]\sum_{j\in\Mfrak_0} \bigl\llVert\hat f_j-f_j^*\bigr\rrVert_n,
\end{eqnarray*}
which, provided $\lambda_n\ge4(C_1+1)r_n$, implies
%
%e30 #&#
\begin{equation}
\label{compatbnd} \sum_{j\in\Mfrak_0^c} \bigl\llVert\hat f_j-f_j^*\bigr\rrVert_n\le3\sum
_{j\in
\Mfrak_0} \bigl\llVert\hat f_j-f_j^*
\bigr\rrVert_n.
\end{equation}
This allows us to apply the compatibility condition,~\ref{as4}(A), to
$\hat f -f^*$. It follows that $s_n^{-1}(\sum_{j\in\Mfrak_0} \llVert
\hat f_j-f^*_j\rrVert _n)^2\le\llVert \hat f-f^*\rrVert
_n^2/\phi^2$, which,
by~(\ref{compatbnd}), yields $s_n^{-1}(\sum_{j=1}^{p_n} \llVert
\hat f_j-f^*_j\rrVert _n)^2\le16\llVert \hat f-f^*\rrVert _n^2/\phi
^2$. Stochastic
bound~(\ref{stochbndlam}) then gives $\sum_{j=1}^{p_n} \llVert \hat f_j-f^*_j\rrVert _n= O(s_n\lambda_n)$, and\vspace*{2pt} hence $\sum_{j=1}^{p_n}
\llVert
\hat f_j-f_{0j}\rrVert _n= O(s_n\lambda_n+s_nd_n^{-2})$, which again
implies the bound in display~(\ref{nonlinthmbnd}). This completes
the proof of Theorem~\ref{nonlinthm}.

Under the assumptions of Theorem~\ref{threshthm}, the error bound in
the statement of Theorem~\ref{nonlinthm} simplifies to $\sum
_{j=1}^{p_n} \llVert \hat f_j-f_{0j}\rrVert _n=O(s_n\lambda_n)$.
Consequently,
on the sets of probability tending to one,
%
%e31 #&#
\begin{equation}
\label{bndprfthrthm} \sum_{j\in\Mfrak_0^c} \llVert\tilde f_j\rrVert_n\le\sum_{j\in\Mfrak
_0^c}
\llVert\hat f_j-f_{0j}\rrVert_n=O(s_n
\lambda_n).
\end{equation}
Using bound~(\ref{bndprfthrthm}) and the fact that $\llVert \tilde f_j\rrVert _n>\lambda_n$ for $j\in\widetilde\Mfrak_n$, we can deduce
$\llvert \Mfrak_0^c\cap\widetilde\Mfrak_n\rrvert =O(s_n)$. This implies
$\llvert \widetilde\Mfrak_n\rrvert \le\llvert \Mfrak_0\rrvert
+\llvert \Mfrak_0^c\cap\widetilde
\Mfrak_n\rrvert =O(s_n)$.
Also note that
\[
\sum_{j\in\Mfrak_0} \llVert\tilde f_j-f_{0j}
\rrVert_n\le\sum_{j\in
\Mfrak_0}\bigl(
\lambda_n+ \llVert\hat f_j-f_{0j}\rrVert
_n\bigr) =O(s_n\lambda_n).
\]
The above bound, together with~(\ref{bndprfthrthm}), yields the
error bound in Theorem~\ref{threshthm}.
\end{longlist}
\end{appendix}

\section*{Acknowledgments}
We would like to thank the Center for Clinical Neurosciences,
University of Texas Health Science Center at Houston for the use of their
MEG data.

\begin{supplement}[id=suppA]
\stitle{Supplementary material for: Functional additive regression\\}
\slink[doi]{10.1214/15-AOS1346SUPP} %[doi,text={...}] - jei reikia
%suskaldyti doi
\sdatatype{.pdf}
\sfilename{AOS1346\_supp.pdf}
\sdescription{Due to space constraints,
the proofs of Theorems~\ref{T2} and~\ref{T3} and Lemma~\ref{lementrcons} are relegated to the
supplement \cite{FAR}.}
\end{supplement}

% imsref loaded by linak, 2015-07-02 13:59:34
%

\printaddresses
\end{document}